\documentclass[12pt,a4paper]{amsart}
\usepackage{graphicx,amssymb}
\input xy
\xyoption{all}
\usepackage[all]{xy}
\usepackage{hyperref}

\newtheorem{thm}{Theorem}[section]
\newtheorem{cor}[thm]{Corollary}
\newtheorem{lem}[thm]{Lemma}
\newtheorem{prop}[thm]{Proposition}
\theoremstyle{definition}

\newtheorem{rem}[thm]{Remark}

\numberwithin{equation}{section}

\newcommand{\ZZ}{\mathbb Z}
\newcommand{\CC}{\mathbb C}
\newcommand{\PP}{\mathbb P}
\newcommand{\FF}{\mathbb F}

\newcommand{\ra}{\rightarrow}

\newcommand{\cA}{\mathcal{A}}

\newcommand{\cJ}{\mathcal{J}}
\newcommand{\cK}{\mathcal{K}}

\newcommand{\tX}{\tilde{X}}
\newcommand{\tY}{\tilde{Y}}
\newcommand{\tZ}{\tilde{Z}}

\newcommand{\tP}{\tilde{P}}
\newcommand{\tT}{\tilde{T}}

\newcommand{\cE}{\mathcal{E}}

\newcommand{\cM}{\mathcal{M}}

\newcommand{\cO}{\mathcal{O}}
\newcommand{\cP}{\mathcal{P}}
\newcommand{\cR}{\mathcal{R}}
\newcommand{\cS}{\mathcal{S}}

\newcommand{\cZ}{\mathcal{Z}}
\newcommand{\cX}{\mathcal{X}}
\newcommand{\cY}{\mathcal{Y}}
\newcommand{\ch}{\mathfrak{h}}

 \DeclareMathOperator{\Ker}{Ker}
  
\DeclareMathOperator{\Pic}{Pic}
 \DeclareMathOperator{\Nm}{{Nm}}

\DeclareMathOperator{\Spec}{Spec}
 \DeclareMathOperator{\Stab}{Stab}

 \DeclareMathOperator{\im}{Im}
 \DeclareMathOperator{\Div}{{Div}}

\begin{document}

\title[ ]{Compactification of the Prym map for non cyclic triple coverings}%
\author{ Herbert Lange,  Angela Ortega}

\address{H. Lange \\ Department Mathematik der Universit\"at Erlangen \\ Germany}
\email{lange@mi.uni-erlangen.de}

\address{A. Ortega \\ Institut f\" ur Mathematik, Humboldt Universit\"at zu Berlin \\ Germany}
\email{ortega@math.hu-berlin.de}

\thanks{We would like to thank the Isaac Newton Institute for Mathematical Sciences, Cambridge. We participated in its program Moduli Spaces 
where most of this paper was written.
The second author was partially supported by Deutsche Forschungsgemeinschaft, SFB 647}

\begin{abstract} 
According to \cite{lo}, the Prym variety of any non-cyclic \'etale triple cover $f: Y \ra X$ of a smooth curve $X$ of genus 2 is a Jacobian variety of dimension 2. This gives a map from
the moduli space of such covers to the moduli space of Jacobian varieties of dimension 2. We extend this map to a proper map $Pr$ of a moduli space 
${}_{S_3}{\widetilde \cM}_2$ of admissible 
$S_3$-covers of genus 7 to the moduli space $\cA_2$ of principally polarized abelian surfaces. The main result is that 
$Pr: {}_{S_3}{\widetilde \cM}_2 \ra \cA_2$ is finite surjective of degree 10.     
\end{abstract}

\maketitle

\section{Introduction}

Let $f:Y \ra X$ denote a non-cyclic cover of degree 3 of a smooth projective curve $X$ of genus 2. The Prym variety $P = P(f)$ of $f$ is by definition the 
complement of the image of the pull-back map of Jacobians $f^*: JX \ra JY$ with respect to the canonical polarization of $JY$. It is easy to see that the canonical 
polarization of $JY$ restricts to the 3-fold of a principal polarization $\Xi$ on $P$. This induces a morphism $Pr$, called Prym map, 
of the moduli space
$\cR_{2,3}^{nc}$ of connected \'etale non-cyclic degree-3  covers of curves of genus 2 into the moduli space $\cA_2$ of principally polarized abelian surfaces. In \cite{lo}
we showed that the image of $Pr$ is contained in the Jacobian locus $\cJ_2$ and moreover that
$$
Pr:  \cR_{2,3}^{nc} \ra \cJ_2
$$
is of degree 10 onto its image. The main aim of this paper is to determine the image of $Pr$, we will see that $Pr$ is not surjective, and  extend this map to a proper surjective map. 

For this it turns out to be convenient to shift the point of view slightly. In \cite[Proposition 4.1]{lo} we saw that taking the Galois closure gives a bijection 
between the set of connected non-cyclic \'etale $f$ covers of above and the set of \'etale Galois covers $h:Z \ra X$, with Galois group the symmetric group $S_3$ of order 6.  
Hence, if we denote by ${}_{S_3}\cM_2$ the moduli space of \'etale Galois covers of smooth curves of genus 2 with Galois group $S_3$ as constructed for example in 
\cite[Theorem 17.2.11]{acg}, we get a morphism which we denote by the same symbol,
$$
Pr:{}_{S_3}\cM_2 \ra \cJ_2,
$$  
and also call the Prym map. Then we use the compactification ${}_{S_3}{\overline \cM}_2$ of ${}_{S_3}\cM_2$ by admissible $S_3$-covers as constructed in 
\cite[Chapter 17]{acg} (based on \cite{acv}) to define the extended Prym map. In fact, consider the following subset of ${}_{S_3}{\overline \cM}_2$:
$$
 {}_{S_3}{\widetilde \cM}_2 := \left\{ [h:Z \ra X] \in {}_{S_3}{\overline \cM}_2 \; \left| \begin{array}{c}
                                                                     p_a(Z) = 7 \; \mbox{and for any node} \; z \in Z \\
                                                                     \mbox{the stabilizer}\; \Stab (z) \; \mbox{is of order 3}

                                                                  \end{array} \right. \right\}.                
$$
Then ${}_{S_3}{\widetilde \cM}_2$ is a non-empty open set of a component of ${}_{S_3}{\overline \cM}_2$ containing the smooth $S_3$-covers ${}_{S_3}\cM_2$.
For any $[h:Z \ra X] \in {}_{S_3}{\widetilde \cM}_2$ let $Y$ denote the quotient of $Z$ by a subgroup of order 2 of $S_3$. We show that the kernel $P = P(f)$
of the map $f:Y \ra X$ is a principally polarized abelian surface. 
Hence we get an extended map $Pr: {}_{S_3}{\widetilde \cM}_2 \ra \cA_2$, which is modular and 
which we 
denote by the same symbol and also call the Prym map. Clearly $P$ does not depend on the choice of the subgroup of order 2. Our main result is the following theorem.\\

\noindent
{\bf Theorem.} {\it The Prym map $Pr: {}_{S_3}{\widetilde \cM}_2 \ra \cA_2$ is a finite surjective morphism of degree $10$.}\\

 In fact, we can be more precise. Consider the following stratification of ${}_{S_3}{\widetilde \cM}_2$:
$$
{}_{S_3}{\widetilde \cM}_2 = {}_{S_3}{\cM}_2 \sqcup R \sqcup S,
$$
where $R$ denotes the set of covers of ${}_{S_3}{\widetilde \cM}_2$ with $X$ singular, but irreducible, and $S$ 
denotes the complement of ${}_{S_3}{\cM}_2 \sqcup R$ in ${}_{S_3}{\widetilde \cM}_2$. As for $\cA_2$,
let $\cE_2$ denote the closed subset of $\cA_2$ consisting of 
products of elliptic curves with canonical principal polarisation. For any smooth curve $C$ of genus 2 and any 3 Weierstrass points $w_1,w_2,w_3$ of $C$ let
$\varphi_{2(w_1+w_2+w_3)}$ denote the map $C \ra \PP^1$ defined by the pencil $(\lambda (2(w_1+w_2+w_3)) + \mu (2(w_4+w_5+w_6))_{(\lambda,\mu) \in \PP^1}$, where $w_4,w_5,w_6$ 
are the complementary Weierstrass points. The map $\varphi_{2(w_1+w_2+w_3)}$ factorizes via the hyperelliptic cover and a $3:1$ map $\bar f: \PP^1 \ra \PP^1$.
With this notation we define the following subsets of $\cJ_2$,
$$
\cJ_2^u:=\{  JC \in \cJ_2 \mid  \  \exists \; w_1,w_2,w_3\mbox{ in } C
\mbox{ such that } \bar{f} \mbox{ is simply ramified} \},
$$
$$
\cJ_2^r:=\{  JC \in \cJ_2 \mid  \  \exists \; w_1,w_2,w_3 \mbox{ in } C
\mbox{ such that } \bar{f} \mbox{ is not simply ramified} \}.
$$
So we have
$$
\cA_2 = \cJ_2^u \cup \cJ_2^r \sqcup \cE_2.
$$ 
We show that the Prym map restricts to finite surjective morphisms $Pr: {}_{S_3}{\cM}_2 \ra \cJ_2^u$, $Pr: S \ra \cE_2$ and to a finite 
morphism $R \ra \cJ_2^r$. We then prove that the extended Prym map is proper and of degree 10, this  implies the theorem. \\

Recall that an even spin curve of genus 2 is a pair consisting of a smooth curve of genus 2 and an even theta characteristic on it and that every curve 
of genus 2 admits exactly 10 even theta characteristics. The degree of the Prym map being 10 suggests that the moduli spaces ${}_{S_3}{\widetilde \cM}_2$ and the 
moduli space of even spin curves should be related. And in fact they are. We will work out details in the forthcoming paper \cite{lo1}. \\

The first part of the paper is devoted to proving that the Prym map $Pr: {}_{S_3}{\widetilde \cM}_2 \ra \cA_2$ is proper. We apply a method used already in \cite{b}
and \cite{f} to show the properness of a Prym map: we consider an open set of ${}_{S_3}{\overline \cM}_2$, which is seemingly bigger than ${}_{S_3}{\widetilde \cM}_2$,
namely the set of $S_3$-covers satisfying condition $(**)$ (see Section 3). In Section 5 we classify these $S_3$-covers and use this to show that this set coincides with 
${}_{S_3}{\widetilde \cM}_2$. In Sections 6 we deduce from this the properness of the extended Prym map. In Sections 7 and 8 we study the restriction of $Pr$
to $R$ and $S$. Finally, Section 9 contains the proof of the above mentioned theorem. \\

We thank Gavril Farkas for some valuable discussions.

\section{Admissible $S_3$-covers}

In this section we recall some notions and results which we need subsequently. For the definitions and results
on admissible coverings we refer to \cite[Chapter 16]{acg} and \cite{acv}.
Let $\cX \ra S$ be a family of connected nodal curves of arithmetic genus $g$ and $d \geq 2$ be an integer. 
A family of {\it degree $d$ admissible covers} of $\cZ$ over $S$ is a finite morphism $\cZ \ra \cX$ such that,
\begin{enumerate}
\item the composition $\cZ \ra S$ is a family of nodal curves;
\item every node of a fiber of $\cZ \ra S$ maps to a node of the corresponding fiber of $\cX \ra S$;
\item away from the nodes $\cZ \ra \cX$ is \'etale of constant degree $d$;
\item if the node $z$ lies over the node $x$ of the corresponding fibre of $\cX \ra S$, 
the two branches near $z$ map to the two branches near $x$ with the same ramification index $r \geq 1$.
\end{enumerate}

If $G$ is a finite group, a $G$-cover $\cZ \ra \cX$ is called a family of {\it admissible $G$-covers} if in
addition to (1) and (2) it satisfies 

\medskip
$\;\;\; (3') \;\;\cZ \ra \cX$ is a principal $G$-bundle away from the nodes;

$\;\, (4') $ if $\xi$ and $\eta$ are local coordinates of the two branches near $z$, any element of 
\hspace*{1.6cm} the stabilizer 
$\Stab_G(z)$ acts as 
$$
(\xi, \eta) \mapsto (\zeta \xi, \zeta^{-1}\eta)
$$
\hspace*{1.6cm} where $\zeta$ is a primitive $r$-th root of the unity for some positive integer $r$. 

\medskip
In the case of $S = \Spec \CC$ we just speak of an admissible degree $d$- (respectively $G$-) cover.
Clearly, $(3')$ and $(4')$ imply $(3)$ and $(4)$. So an admissible $G$-cover is an admissible $d$-cover with $d = |G|$.
In the case of an admissible $G$-cover,
the ramification index at any node $z$ over $x$ equals the order of the stabilizer of $z$ and depends only on $x$. 
It is called the {\it index} of the $G$-cover $\cZ \ra \cX$ at $x$.
Note that, for any admissible $G$-covering $Z \ra X$, the curve $Z$ is stable if and only if and only if $X$ is stable. 

In this paper we are interested in the case $G = S_3$, with
$$
S_3 := \langle \sigma, \tau \;|\; \sigma^3= \tau^2 = \tau \sigma \tau \sigma = 1 \rangle.
$$ 
If $h: Z \ra X$ is any $S_3$-covering of nodal curves, then all curves in the following diagram
\begin{equation} \label{diag2.1}
\xymatrix{
 & Z \ar[dl]_{p} \ar[dr]^{q} \ar[dd]^{h} & \\
Y = Z/\langle \tau \rangle   \ar[dr]_{f}  & & D = Z/\langle \sigma \rangle  \ar[dl]^{g} \\
 & X = Z/S_3  & \\
    }
\end{equation}
have only ordinary nodes too. Here $p$ and $g$ are of degree 2, $q$ is cyclic of degree 3 and $f$ 
is non-cyclic of degree 3. There is, of course, an analogous diagram for any family of $S_3$-coverings.

Let $z$ be a node of $Z$ such that every element in $\Stab z := Stab_{S_3}(z)$ does not exchange the 2 branches of $Z$ at $z$,
then the subgroup $\Stab z$ is cyclic (\cite[p. 529]{acg}). This follows from the fact that $\Stab z$ injects into the automorphism group 
of the tangent space to a branch at $z$. It is then easy to see that any node $z$ of $Z$ is of one of the following types:
\medskip
\begin{enumerate}
\item $|\Stab z| = 1$. The orbit $S_3(z)$ consists of 6 nodes and its image $x = h(z)$ is a node of $X$.
\item $|\Stab z| = 2$ and the generator of $\Stab z$ does not exchange the 2 branches of $Z$ at $z$. 
      The orbit $S_3(z)$ consists of 3 nodes and its image $x = h(z)$ is a node of $X$.
\item $|\Stab z| = 3$. The generator of $\sigma$ of $\Stab z$ acts on each branch of $Z$ in $z$,
      the orbit $S_3(z)$ consists of 2 nodes and its image $x = h(z)$ is a node of $X$.  
\item $|\Stab z| = 2$ and the generator of $\Stab z$ does exchanges the 2 branches of $Z$ at $z$. 
      The orbit $S_3(z)$ consists of 3 nodes and its image $x = h(z)$ is smooth in $X$.
\item $\Stab z = S_3$. In this case $\tau$ exchanges the 2 branches of $Z$ at $z$, the orbit $S_3(z)$ consists of $z$ alone and 
$x = h(z)$ is smooth in $X$. If  $\tau$ does not exchanges the branches\footnote{Note that  if $\tau$ exchanges the branches, so do all 
the other elements of order 2.}  of $Z$ at $z$, then $\Stab z$ is also the stabilizer of $z$ on a branch, but the stabilizers on a smooth 
curve are always cyclic, so we can exclude this case.

\end{enumerate}
We call the nodes {\it of type} (1), \dots, (5) respectively. 

Note that if the $S_3$-covering $Z \ra X$ is admissible, then 
every node of $Z$ is of type (1), (2) or (3), because when the node is of type (4) or (5) the map $h$ does not verify condition (2) 
of the definition. If $z$ is of type (1), then all maps in \eqref{diag2.1} are \'etale near $z$ and its images. 
If $z$ is a node of type (2), then $q$ and $f$ are \'etale 
near $z$ and $p(z)$ and $p$ and $g$ are ramified at both branches near $z$ and $q(z)$. 
If $z$ is a node of type (3), then $p$ and $g$ are \'etale near $z$ and $q(z)$ and  
$q$ and $f$ are ramified of index 3 at both branches near $z$ and $p(z)$. 

In order to describe the norm map of $f:Y \ra X$ (note that the curves $X$ and $Y$ are not necessarily irreducible) 
we need the following description of the divisors. We have 
$$
\Div (Y) = \bigoplus_{x \in Y_{sm} } \ZZ \cdot x + \bigoplus_{y \in Y_{sing}} \cK_{Y,y}^* / \cO_{Y,y}^*
$$
and similarly for $\Div (X)$,
where $\cK_Y$ (resp. $\cK_X$) is the ring of rational function of $Y$ (resp. $X$). 
Let $n_Y: \tY \ra Y$ and $n_X: \tX \ra X$ be the normalizations. For a node $y$ of $Y$ denote 
$n_Y^{-1}(y) = \{ y_1, y_2 \}$ and $\nu_i$ the valuation of $y_i$ for $i = 1,2$. 
We obtain an isomorphism (see \cite[\S 3]{b})
$$
\cK_{Y,y}^* / \cO_{Y,y}^* \simeq \CC^* \times \ZZ \times \ZZ \quad \mbox{defined by} \quad
\phi \mapsto \left( \frac{\phi(y_1)}{ \phi(y_2)}, \nu_1(\phi), \nu_2(\phi) \right)
$$ 
and similarly for $\cK_{X,f(y)}^* / \cO_{X,f(y)}^*$. With these identifications we have the following lemma. 
\begin{lem} \label{pushfwd} Let $y$ be a node of $Y$ and $(\gamma, m,n) \in \cK_{Y,y}^* / \cO_{Y,y}^*$.

{\em (a)} If $y$ is of type $(1)$ or $(2)$, then $f_*(\gamma, m,n) = (\gamma,m,n)$; 

{\em (b)} If $y$ is of type $(3)$, then $f_*(\gamma, m,n) = (\gamma^3,m,n).$
\end{lem}

\begin{proof}
(a) is a consequence of the fact that $f$ is \'etale near $y$.
(b) follows from diagram \eqref{diag2.1} using the facts that the maps $p$ and $g$ are \'etale near the corresponding nodes and  
the analogous statement for the cyclic map $q$ which was shown in \cite[p.64]{f}.
\end{proof}

Let $JX$ and $JY$ denote the (generalized) Jacobians of $X$ and $Y$. The closed points in the Jacobian $JX$ (resp. $JY$)  
can be identified with the isomorphism classes of line bundles on $X$  (resp. on $Y$) of multidegree $(0, \ldots , 0)$.
Let $\Nm_{f}: \Pic(Y) \ra \Pic(X)$ be the norm map. Since the norm of $\cK_Y /  \cK_X$ maps $\cO_Y$ into $\cO_X$ we get the 
diagram of exact sequences
\begin{equation} \label{}
\xymatrix{
	      \cK_Y^*  \ar[d]^{N_{\cK_Y/\cK_X}} \ar[r] &  \Div (Y) \ar[d]^{f_*} \ar[r] & \Pic(Y) \ar[r] \ar[d]^{\Nm_f} & 0  \\
               \cK_X^*   \ar[r] &  \Div (X)  \ar[r] & \Pic(X) \ar[r] & 0 
    }
\end{equation}
The norm map induces a morphism $\Nm_f: JY \ra JX $. We define the Prym variety associated to $f$ as 
$$
P(f):= (\Ker \Nm_{f})^0.
$$ 
In general the kernel of the norm can have several conected components. However, consider the following condition
\begin{equation*}
(*) \left\{ \begin{array}{l}
      \mbox{let $h:Z \ra X$ be an admissible $S_3$-cover}\\
       \mbox{such that all nodes of} \; Z \; \mbox{are of type (3)}.
\end{array}        
\right.
\end{equation*}
Then we have,

\begin{lem} \label{lemma2.2}
Let $h: Z \ra X$ be an $S_3$-cover satisfying $(*)$. Then 
$\Ker \Nm_{f}$ is an abelian subvariety $P$ of $JY$.
\end{lem}     

\begin{proof}
Let $n_3$ be the number of nodes and $s$ be the number of irreducible components of $X$. Then we have the following exact
diagram. 
\begin{equation} \label{diagram}
\xymatrix{
0 \ar[r] & T_3 \ar[r] \ar@{^{(}->}[d] & K \ar[r] \ar@{^{(}->}[d] & R \ar[r] \ar@{^{(}->}[d] & 0\\
0 \ar[r] & T_Y \ar@{->>}[d]^{Nm_{f}} \ar[r] & JY \ar@{->>}[d]^{\Nm_{f}} \ar[r]^{n_Y^*} & J\tY \ar@{->>}[d]^{\Nm_{\tilde{f}}} \ar[r] & 0\\
0 \ar[r] & T_X \ar[r]  & JX \ar[r]^{n_X^*}  & J\tX \ar[r] & 0\\
   }
\end{equation}
with
$$
T_Y \simeq T_X \simeq {\CC^*}^{n_3 -s +1},
$$
(for the last isomorphism see \cite[1.1.3]{c}).
The kernel $R$ of $\Nm_{\tilde{f}}$ is an abelian subvariety of $J \tilde Y$, since $\tilde Y$ and $\tilde X$ are 
disjoint unions of smooth projective curves and on every component $Nm_{\tilde{f}}$ is non-cyclic (and ramified) of degree 3.
Moreover, according to Lemma \ref{pushfwd}  ${\Nm_{f}}_{|_{T_Y}}$ is surjective and the kernel $T_3$ of $\Nm_{f}$ is isomorphic to 
$(\ZZ/3 \ZZ)^{n_3 - s +1}$.  This implies that the kernel $K$ of $\Nm_{f}$ is a compact subgroup of $JY$.

It remains to show that $K$ is connected, i.e. $K = P$. In order to show this we consider a connected family of $S_3$-covers 
whose general member is smooth and which contains the given cover $f:Z \ra X$ as a special fiber. 
We get a family of maps $\Nm_{f}$ in an obvious way and thus a
family $\pi: \mathcal{F} \ra B$ whose fibres are the kernels $\Ker\Nm{f}$. We may assume the base $B$ is smooth. 
Applying the Stein factorization to $\pi$ we get a morphism $\pi': \mathcal{F} \ra B'$ with connected fibres and a finite
morphism $b: B' \ra B$, such that $\pi = b \circ \pi'$. By assumption, the generic fibre of $\pi$ is connected since $\Ker\Nm{f}$ is 
connected for smooth coverings, hence $b$ must be the identity and $\pi= \pi'$. Therefore, the special fibre $K$ is also connected. 

\end{proof}

Let $L \in \Pic^3(Y)$ be a fix line bundle on $Y$.  According to \cite{b} (see also \cite{c}) the set 
$\Theta_L :=  \{ M\in JY \mid \ H^0(Y, L\otimes M) \geq 1 \}$
defines a theta divisor on $JY$, which is linearly equivalent to $(n_Y^*)^{-1}(\Theta_{L'}) $, where $\Theta_{L'}$ is a theta 
divisor on $J\tY$. 

\begin{prop}
Let $h: Z \ra X$ be an $S_3$-cover satisfying $(*)$ and let $\rho : P \ra \hat{P}$ be the polarization 
induced by the principal polarization of $JY$ . 
Then 
$$
|\Ker \rho  \ | = 3^{2p_a(X)}.
$$
\end{prop}
\begin{proof}
 Consider the isogeny 
$$
\beta: P \times J\tX \ra J\tY, \; (L,M) \mapsto n_Y^*L \otimes \tilde{f}^*M.
$$ 
Then $\Ker(\beta )$ is a maximal isotropic subgroup of the kernel of the polarization on $P\times J\tX$ 
given by the pullback of the principal polarization on  $J\tY$. Consider again the commutative diagram \eqref{diagram}
with $K = P$ according to Lemma \ref{lemma2.2}. 
 An element $(L, M) \in P
\times J\tX$ is in the kernel of $\beta$ if and only if 
$$
n_Y^*L \otimes \tilde{f}^*M \simeq \cO_{\tY}.
$$ 
Let $M= n_X^*M' $, for some $M' \in JX$,
then $L \otimes f^*M' \in \Ker n_Y^* = T_Y $. 
Since $T_3$ is finite, we have $T_Y=f^*(T_X)$. This implies that $L = f^*M''$ for some $M'' \in JX$. 
 Since $\Nm_{f} (L) \simeq \cO_X$ we obtain  $(M'' )^3 \simeq \cO_X $,  which implies 
 $$
 L^3 \simeq \cO_Y \quad \mbox{and} \quad M^3 \simeq \Nm_{\tilde{f}} \tilde{f}^*M \simeq \Nm_{\tilde{f}} n_Y^* L^{-1} 
 = n_X^* \Nm_{f} L^{-1} \simeq \cO_{\tX}.
 $$
So $\Ker \beta \subset P[3] \times J\tX[3] $ and therefore $\Ker \rho \subset P[3]$.  Moreover,
\begin{eqnarray*}
 \tilde{f}^*M  & = & n_Y^*L^{-1} \\
  & = & n_Y^* f^*( M'' )^{-1}  \\
    & = & \tilde{f}^* n_X^* ( M'' )^{-1} .
\end{eqnarray*}
Since $\tilde{f}$ is non-cyclic on each irreducible component of $\tY$,  $\tilde{f}^*$ is injective. Then 
$M  =  n_X^* ( M'' )^{-1} $. This shows that
 $$
 \Ker \beta =\{ ( f^*a, n_X^* a^{-1})  \mid a \in JX[3] ,  \  f^*a \in P \}.
 $$
Since $f^*$ is injective when $X$ is singular or $n_X^*$ is injective when $X$ is non-singular, we conclude that 
$\Ker \beta \simeq \{ a \in JX[3] \mid f^*a \in P \}.$ 
Since $\Nm \circ f^*$ is the multiplication by 3,  $f^*(JX[3]) \subset \Ker \Nm_{f} = P $, by 
the previous lemma. Hence we obtain
$$ 
\Ker \beta \simeq JX[3].
$$ 
Let $t = \dim T_Y =\dim T_X = n_3 -s +1$ and denote $\mu: J\tY \ra  \widehat{J\tY}$ the canonical principal polarization. 
Then the dimension over $\FF_3$ of the kernel of the pullback of $\mu $ on $P \times J\tX$ is
 \begin{eqnarray*}
  \dim_{\FF_3} \Ker (\hat{\beta} \circ \mu \circ \beta) & = &  2 \dim_{\FF_3} \Ker \beta \\
    & = &  2 \dim_{\FF_3} JX[3]\\
    &=&  2( 2p_a(X) - t ) 
 \end{eqnarray*}
 On the other hand we clearly have
 $$
 \dim_{\FF_3} \Ker (\hat{\beta} \circ \mu \circ \beta)  =   \dim_{\FF_3} (\Ker \rho \times J\tX[3] ).
 $$
 Since $\dim_{\FF_3} J\tX[3] = 2(g(\tX)) = 2(p_a(X) -t)$ we obtain  $ \dim_{\FF_3} \Ker \rho= 4p_a(X) -2t - \dim_{\FF_3} J\tX[3] =  2p_a(X) . $
\end{proof}

\begin{cor}  \label{cor2.4}
Let $h: Z \ra X$ be an $S_3$-cover satisfying $(*)$.  Then the polarization $\rho : P \ra \hat{P}$ is three times a principal polarization
if and only $\dim P = p_a(X)$.
\end{cor} 
\begin{proof}
We have  $\Ker \rho $ is a subset of $ P[3]$ of cardinality $3^{2 p_a(X)}$ and this is equal to the cardinality 
of $P[3]$ if and only if $\dim P = p_a(X)$. 
\end{proof}

\begin{cor} \label{cor2.5}
Let $h: Z \ra X$ be an $S_3$-cover satisfying $(*)$ and let $X_1, \dots, X_s$ be the irreducible components of $X$. Then
the canonical polarization of $P$ is three times a principal polarization if and only if 
$$
\sum_{i=1}^s p_g(X_i) = s - n_3 + 1.
$$
where as above $n_3$ denotes the number of nodes of $X$. In this case $\dim P = 2$.
\end{cor}

\begin{proof}
Since $s$ is also the number of components of $Y$ and $n_3$ the number of nodes of $Y$, we have
$$
p_a(Y) = \sum_{i=1}^s p_g(Y_i) + n_3 -s + 1
$$
and 
$$
p_a(X) = \sum_{i=1}^s p_g(X_i) + n_3 -s + 1.
$$
Hence $\dim P = p_a(Y) - p_a(X) = p_a(X)$ if and only if
\begin{equation} \label{eq2.2}
\sum_{i=1}^s p_g(Y_i) = 2 \sum_{i=1}^s p_g(X_i) + n_3 -s + 1.
\end{equation}
The covering $Y\ra X$ is doubly ramified exactly over each branch of the nodes of $X$. So if $2r_i$ denotes 
the order of the ramification divisor of the normalization ${\tilde Y}_i \ra {\tilde X}_i$ of $Y_i \ra X_i$ we have 
$$
\sum_{i=1}^s r_i = 2n_3.
$$
On the other hand, by the Hurwitz formula, $p_g(Y_i) = 3p_g(X_i) -2 + r_i$. So \eqref{eq2.2} gives
$$
3\sum_{i=1}^s p_g(X_i) -2s + 2n_3 = 2\sum_{i=1}^s p_g(X_i) +n_3 -s + 1.
$$
This implies the first assertion. Finally we have
$$
\dim P = p_a(X) = \sum_{i=1}^s p_g(X_i) + n_3 -s + 1 = s- n_3 + 1 +n_3 -s +1 = 2.
$$
\end{proof}

\section{The condition $(**)$}

As in \cite[Section 5]{b} we shall need to study the Prym variety $P$ of $\pi$ under a more general assumption than 
the hypothesis $(*)$ of the previous section. Let $Z$ be a connected curve with only ordinary nodes and an $S_3$-action
and consider the diagram \eqref{diag2.1}. Throughout we assume condition $(4')$ of the definition of an $S_3$-cover, i.e.
that at each node of type (3), if one local parameter is 
multiplied under $\sigma$ by $\zeta$, then the other is multiplied by $\zeta^2$ for some third root of unity $\zeta$ (see \cite[Remark 2.1]{f}).  
Moreover, we assume that the action satisfies the following condition
\begin{equation*}
(**) \left\{ \begin{array}{l}
        P \; \mbox{is an abelian variety};\\
        \sigma \; \mbox{and} \;\tau \; \mbox{are not the identity on any component of $Z$};\\
        p_a(Z) = 6p_a(X) -5.
\end{array}        
\right.
\end{equation*}
The number
$$
n_i : = |\{\mbox{nodes} \; z \; \mbox{of} \; Z \; \mbox{of type (i)} \}| \cdot \frac{|\Stab z|}{6}
$$
obviously is a non-negative integer. For $i = 1,2,3$ it coincides with the number of nodes of $X$ of type (i). We set
$$
c_i : = |\{ \mbox{components} \; Z_i \; \mbox{of} \; Z\; \mbox{with}\;  |\Stab Z_i| = i \}| \cdot \frac{i}{6} 
\quad \mbox{for} \; i = 1,2,3,6
$$
so the number of components of $X$ is $c_1+c_2+c_3+c_6$. Finally define
$$
r_i := |\{ \mbox{smooth points}\; z \; \mbox{of} \; Z \; \mbox{with} \; |\Stab z| = i \}| \quad \mbox{for} \; i = 2,3.  
$$

\begin{prop} \label{lem3.1}
The assumptions $(**)$ are equivalent to 
$$
\left\{
                                             \begin{array}{l}
                                             r_2= r_3 = n_4 = n_5 = 0,\\
                                             2n_1 +n_2 = 2c_1 + c_2.\\
                                             \end{array}                                \right.
$$                                                                                        
In particular $h:Z \ra X$ is an admissible $S_3$-cover.                                           
\end{prop}

\begin{proof}
Let $\tZ$ and $\tX$ denote the normalizations of $Z$ and $X$. The induced covering $\tilde{h}:\tZ \ra \tX$ is ramified
exactly at points of $\tZ$ lying over non-singular points fixed by $\sigma$ or an element of order 2 and at nodes of 
types (2), (3) and (5).
Hence, by the Hurwitz formula:
$$
p_a(\tZ) = 6 p_a(\tX) -5 + \frac{r_2}{2} + r_3 + 3n_2 + 4n_3 + 2n_5.
$$
So
\begin{eqnarray*}
p_a(Z) &=& p_a(\tZ) + 6n_1 + 3n_2 + 2n_3 + 3n_4 + n_5\\
& = &  6p_a(\tX) - 5 + \frac{r_2}{2} + r_3 + 6n_1 + 6n_2 + 6n_3 + 3n_4 + 3n_5.
\end{eqnarray*}
The nodes of $X$ come from nodes of $Z$ of types (1), (2) and (3), hence:
$$
p_a(X) = p_a(\tX) + n_1 + n_2 + n_3.
$$ 
Therefore
$$
p_a(Z) = 6p_a(X) -5 + \frac{r_2}{2} + r_3 + 3n_4 + 3n_5.
$$
Hence the condition $p_a(Z) = 6p_a(X) -5$ is equivalent to $r_2 = r_3 = n_4 = n_5 = 0$.

In order to express the condition that $P$ is an abelian variety (in particular $P$ is connected), let $\tY$ denote the normalization 
of $Y$.  Consider the commutative diagram
\begin{equation*} \label{}
\xymatrix{
0 \ar[r] & \tT \ar[d] \ar[r] & JY \ar[d]^{N_{f}} \ar[r] & J\tY \ar[d]^{N_{\tilde{f}}} \ar[r] & 0\\
0 \ar[r] & T \ar[r] & JX \ar[r] & J\tX \ar[r] & 0.\\
    }
\end{equation*}
Note that, since the norm maps $\tT \ra T$ and $J\tY \ra J\tX$ are surjective, $N_f$ is surjective. From the commutative
of the diagram above follows that $P$ is an abelian variety if and only if $\dim \tT = \dim T$.
Now we have
$$
\dim \tT = 3n_1 + 2n_2 + n_3 -3c_1 - 2c_2 -c_3 -c_6 + 1.
$$
For the summand $2n_1$ (and similarly $2c_1$) note that for a transitive action of $S_3$ on a set of 3 elements, $\tau$ 
fixes one element and exchanges the 2 other elements (see Lemma \ref{3points}). Similarly we have 
$$
\dim T = n_1 + n_2 + n_3 - c_1 - c_2 - c_3 - c_6 + 1.
$$
So $\dim \tT = \dim T$ if and only if $2n_1 +n_2 = 2c_1 + c_2$. 
\end{proof}

Apart from the curves of Theorem 2.3 and its corollaries there are some other curves with $S_3$-action which lead to 
principally polarized Prym varieties. These are the analogues of those occurring in \cite[Section 5]{b} and \cite[Section 2]{f}. 
However, in the case of interest for us, i.e. $X$ of genus 2, they do not occur (see Corollary \ref{cor5.7} below).

\section{Some auxiliary results}

Let $X$ be a stable curve of arithmetic genus 2. In the next section we determine the $S_3$-covers $h:Y \ra X$ satisfying condition $(**)$. 
For this we need some lemmata which we collect in this section.

\begin{prop} \label{autgenus3}
Every smooth non-hyperelliptic curve $Z$ of genus $3$ with $S_3$-action has quotient $Z/S_3 \simeq \PP^1$. The automorphism $\sigma$ has $2$ and $\tau$ has $8$
fixed points on $Z$.
\end{prop}

\begin{proof}
Every non-hyperelliptic curve $Z$ of genus 3 with $S_3$-action has an equation (see \cite{v})
$$
z_0^3z_2 + z_1^3z_2 +z_0^2z_1^2 + a z_0z_1z_2^2 + bz_2^4 = 0.
$$
The group $S_3$ is generated by
$$
\sigma: \left\{ \begin{array}{lll}
                z_0 & \mapsto & \zeta z_0\\
                z_1 & \mapsto & \zeta^2 z_1 \\
                z_2 & \mapsto & z_2
                \end{array} \right.
\quad \mbox{and} \quad 
\tau: \left\{ \begin{array}{lll}
              z_0 & \mapsto & z_1 \\
              z_2 & \mapsto & z_2
              \end{array} \right.
$$
where $\zeta$ denotes a primitive third root of unity.
The quotient $D = Z/\langle \sigma \rangle$ is an elliptic curve, since $\sigma$ has exactly 2 fixed points, namely $(1:0:0)$ and 
$(0:1:0)$. An equation of $D$ is 
$$
x_0^2x_1x_2 + x_0x_1^2x_2 + x_0^2x_1^2 + ax_0x_1x_2^2 + bx_2^4 = 0
$$
and the map $Z \ra D$ is given by $x_0 = z_0^3,\; x_1 = z_1^3, \; x_2 = z_1z_2z_3$. The involution $\tau$ 
induces an involution $\overline{\tau}$ on $D$, which is given by $x_0 \mapsto x_1, \; x_2 \mapsto x_2$.
Since $\overline{\tau}$ admits exactly 4 fixed points with multiplicities, the quotient $Z/S_3 = D/ \langle \overline{\tau} \rangle $
is of genus 0.             
\end{proof}

\begin{prop} \label{authypgen3}
 Every smooth hyperelliptic curve $Z$ of genus $3$ with $S_3$-action has an elliptic curve as quotient $Z/S_3$. 
There is a one-dimensional family of such curves $Z$.
\end{prop}

\begin{proof}
Every hyperelliptic curve $Z$ of genus 3 with $S_3$-action has an affine equation (see \cite{i})
$$
y^2 = x(x^3 - 1)(x^3 - a^3)
$$
with $a^3 \neq 0,1$. The group $S_3$ is generated by
$$
\sigma: \left\{ \begin{array}{lll}
                x & \mapsto &\zeta x\\
                y & \mapsto & \zeta^2y 
                \end{array} \right.
\quad \mbox{and} \quad 
\tau: \left\{ \begin{array}{lll}
              x & \mapsto & ax^{-1} \\
              y & \mapsto & -a^2x^4y
              \end{array} \right.
$$
The fixed points of $\sigma$ are a subset of the set of branch points $B =\{0,\infty,1,\zeta,\zeta^2,a,a\zeta,a\zeta^2 \}$.
Hence the quotient map $Z \ra E = Z/\langle \sigma \rangle $ is ramified exactly at the 2 points over $0$ and $\infty$, which
implies that $E$ is an elliptic curve. Since $\tau$ permutes $0$ and $\infty$ as well as the 2 orbits $\{1,\zeta,\zeta^2\}$
and $\{a,a\zeta,a\zeta^2\}$ and $\langle \sigma\rangle $ is a normal subgroup of $S_3$, $\tau$ induces a fixed point free 
automorphism $\overline{\tau}$ of $E$. Hence $Z/S_3 = E/\langle \tau\rangle $ is an elliptic curve.
\end{proof}

\begin{prop} \label{quot}
 There is no smooth genus $2$ curve with automorphism $\sigma$ of order $3$ such that the quotient $Z / \langle \sigma \rangle$
is an elliptic curve.
\end{prop}

\begin{proof}
According to \cite{bo} every smooth curve of genus 2 admitting an automorphism of order 3 has an affine equation
$$
y^2 = (x^3 - a^3)(x^3 - a^{-3})
$$
with $a \neq 0$ and not a $6$th root of unity and $\sigma$ is given by
$$
                x \mapsto \zeta x \quad \mbox{and} \quad
                y \mapsto  y. 
$$
  with a primitive third root of unity $\zeta$.
Hence the quotient $Z / \langle \sigma \rangle$ satisfies the equation $y^2 = (v - a^2)(v - a^{-1})$ which implies that it is rational.  
\end{proof}



\begin{lem} \label{smpoint}
Let $h: Z \ra X$ be an $S_3$-cover of connected nodal curves over $X= X^1 \cup X^2$, where $X_1$ and
$X_2$ are the irreducible components. If $x \in X^1 \cap X^2$,  then
$Z^i := h^{-1}(X^i)$ is smooth at each point of $h^{-1}(x)$, for $i=1,2$.
\end{lem}

\begin{proof}
Let $x \in X^1 \cap X^2$. Then $h^{-1}(x) \subset Z^1 \cap Z^2$ and each point on the fibre of $x$ is a node such that one branch
belongs to a component of $Z^1$and the other to a component of $Z^2$. Thus $Z^i := h^{-1}(X^i)$ is smooth at each point of 
$h^{-1}(x)$.
\end{proof}

\begin{lem} \label{3points}
Let $S$ denote a set consisting of $3$ points $Z_1,Z_2,Z_3$, with a nontrivial $S_3$-action. Then the points may be labelled in  such a 
way that $\tau$ fixes $Z_1$ and exchanges $Z_2$ and $Z_3$ and either
\begin{enumerate}
\item $\sigma(Z_i) = Z_{i+1}$ for $i = 1,2,3$ (where $Z_4 = Z_1$) or 
\item  $\sigma(Z_i) = Z_i$ for $i=1,2,3$.
\end{enumerate}
\end{lem}

The proof is straightforward and will be omitted. 

\begin{cor} \label{3comp}
 Let $Z = Z_1 \cup Z_2 \cup Z_3$ be a nodal curve such that $Z_i \cap Z_{i+1} \; (Z_4 = Y_1)$ consists of one point $n_i$ for $i = 1,2,3$.
Then any $S_3$-cover $h:Z \ra X$ does not satisfy condition $(**)$.
 \end{cor}

\begin{proof} 
Consider the induced $S_3$-action on the dual graph of $Z$ and let the notation be as in the previous lemma.

 First assume case $(1)$, i.e. $\sigma$ permutes the $Z_i$ and $\tau (Z_1)= Z_1$. Then
$$
\tau(n_2) = \tau (Z_2 \cap Z_3 )= Z_3 \cap Z_2 = n_2.
$$ 
Since $\tau(Z_2) = Z_3$, the action on $n_2$ is of type (4), that is, $n_4\neq 0$, contradicting Proposition 3.1. 

Assume case (2), so $\sigma(Z_i)= Z_i$ for $i=1,2,3$. In this case the node $n_2 = \tau (n_2)$ is of type (5), contradicting Proposition 
3.1. 
\end{proof}

\begin{lem}
Let $\Gamma$ be a connected graph with a transitive $S_3$-action consisting of $6$ vertices and $6$ nodes. Then the vertices $Z_i$ and the nodes $n_i$
can be labeled in such a way that 
$$
\sigma = (Z_1,Z_3,Z_5)(Z_2,Z_4,Z_6) \quad \mbox{and} \quad \tau = (Z_1,Z_2)(Z_3,Z_6)(Z_4,Z_5)
$$
on the vertices, and $n_i =  Z_i \cap Z_{i+1}$ where $Z_7 = Z_1$, with
$$
\sigma = (n_1,n_3,n_5)(n_2,n_4,n_6) \quad \mbox{and} \quad \tau = (n_1)(n_4)(n_2,n_6)(n_3,n_5).
$$
Here the notation is the cycle-notation of permutations.

In particular, up to isomorphisms, there is only one transitive $S_3$-action on $\Gamma$.
\end{lem}

The proof is straightforward and will be omitted.

\begin{cor} \label{6comp}
Let $Z$ be a connected nodal curve with an $S_3$-action consisting of $s=6$ irreducible components and $\delta \leq 6$ nodes such that the quotient $X = Z/S_3$ is irreducible.
Then the image of any node of $Z$ is a smooth point of $X$. In particular, $h: Z \ra X$ is not an admissible $S_3$-covering.
\end{cor}
\begin{proof}
The $S_3$-action induces a transitive action on the dual graph $\Gamma$ of $Z$. Since there is no connected graph with 6 vertices and at most 5
edges, we necessarily have $\delta = 6$. Let the notation be as in the previous lemma with components $Z_i$ and nodes $n_i$.
The quotient $D:= Z / \langle \sigma\rangle$ consists of 2 components $D_1= q(Z_1)= q(Z_3)= q(Z_5)$ and 
$D_2=q(Z_2)= q(Z_4)= q(Z_6)$ intersecting transversally in 2 points $q(n_1)$ and $q(n_2)$. The involution $\bar{\tau}$, induced
by $\tau$ in $D$, interchanges $D_1$ and $D_2$ so the quotient $X= D / \bar{\tau} $ is smooth.
\end{proof}

\begin{prop} \label{6components}
Let $h: Z \ra X$ be an $S_3$-cover of connected nodal curves, with $X$ stable of arithmetic genus $2$ consisting of $2$ components  
$X^1$ and $X^2$ intersecting in one point.

If the cover satifies condition $(**)$, then $Z^j = h^{-1}(X^j)$ consists of at most $3$ irreducible components for $j=1$ and $2$.
\end{prop}

\begin{proof}
Suppose $Z^1$ has 6 irreducible components. If the curves $X^1$ and $X^2$ are smooth, the only nodes of $Z$ are on the fibre
over the node of $X$, then the 6 components 
are disjoint according to Lemma \ref{smpoint}. In order to have $Z$ connected, $Z^2$ has to be irreducible and smooth. 
Hence $Z^2 \ra X^2$ is an \'etale map, so $g(Z^2)=1$. 
This is a contradiction, since there is no elliptic curve with a non-trivial $S_3$-action.
 
The components $X^1$ and $X^2$ cannot have more than 1 node, since $X$ is stable of arithmetic genus 2. 
If one of the components of $X$, say $X^1$, has an ordinary double point $x_1$, then $Z^1$ has 6 ordinary double points,  because
the maps $Z^1_i \ra X^1$ are birational. But then the cover $Z^1 \ra X^1$ satisfies the conditions of Corollary \ref{6comp}.
So the image of the node is a smooth point of $X$ contradicting the admissiblity of the covering $h$.
\end{proof}

\section{$S_3$-covers satisfying $(**)$}

In this section we determine the $S_3$-covers $Z \ra X$ satisfying condition $(**)$ with $p_a(X) = 2$. There are 6 types of non-smooth stable curves 
of arithmetic genus 2 which will be considered separately. In the first 2 propositions we assume that $X$ is irreducible.
So let $h:Z \ra X$ be an $S_3$-cover satisfying $(**)$ and 
denote 
$$
Z = \cup_{i=1}^s Z_i
$$ 
with irreducible components $Z_i$. Since $S_3$ acts transitively on the set of components,  
the number $s$ can take the values $1,2,3$ or 6. Moreover, we have the following formula
\begin{equation} \label{eq5.1}
7 = p_a(Z) = \sum_{i=1}^s g_i -s + \delta +1 = s(g_1-1) + \delta + 1,
\end{equation}
where $g_i$ denotes the geometric genus of $Z_i$ and $\delta$ the number of nodes of $Z$. Note that $g_i = g_1 \geq 1$ for all $i$.

Suppose first that $X$ is of geometric genus 1 with one node $x$. 
If $r$ denotes its ramification index,
we have as usual
$$
r \delta = 6.
$$

\begin{prop} \label{prop5.1}
Let $h: Z = \bigcup_{i=1}^s Z_i \ra X$ be an $S_3$-cover of an irreducible curve $X$ of geometric genus $1$ with one 
ordinary double point. 
Then only in the following cases the cover satisfies condition $(**)$:

{\em (a)} $s=1, \; r = 3, \; \delta = 2$. $Z$ is an irreducible curve of geometric genus $5$ with $2$ nodes.

{\em (b)} $s=2, \; r=3, \; \delta = 2$. The normalization of $Z$ consists of $2$ copies of a hyperelliptic curve of genus $3$, admitting an 
automorphism $\sigma$ of order $3$ with $2$ fixed points which are glued together transversally at opposite fixed points of $\sigma$.

In both cases there is a $2$-dimensional family of such coverings. If $\tilde f: \tilde Y \ra \tilde X$ denotes the normalization of $f$,
 in both cases
the Prym variety $P$ is an extension of the Prym variety of $\tilde f$ by the group $\ZZ/3\ZZ$.
\end{prop}

\begin{proof} Let $f$ be an $S_3$-covering satisfying $(**)$. Then $s \neq 3$ according to  Corollary \ref{3comp} 
and  $s \neq 6$ by Corollary \ref{6comp}.
So suppose $s \leq 2$. Then $c_1 = c_2 = 0$ and hence $n_1 = n_2 =0$. This implies $r =3$ and $\delta =2$. If $s =1$, we are in case (1)
by the Hurwitz formula and if $s =2$ we are in case (b) again using the Hurwitz formula.

Concerning the existence statement, there is a 2-dimensional family of curves $X$ of genus 1 with one node. Hence it suffices to show
that, for a fixed such curve $X$, there exist only finitely many admissible $S_3$-covers in the cases (a) and (b).

As for case (a), according to \cite[§3,1,(1)]{kk}, every elliptic curve $E$ admits finitely many  $S_3$-covers of genus 5. The ramification divisor 
of the normalization $\tilde h: \tilde Z \ra \tilde X$ of $h$ is of degree 8. Since 8 is not divisible by 3, $\tilde h$ is doubly ramified at 2 pairs of points.
Gluing them pairwise together we get case (a). 

In case (b),
there are finitely many cyclic coverings $Z^1$ of the normalization of $X$ of degree and genus 3, ramified exactly over the 
2 preimages of the node. For each such cover $Z^1$ take a copy $Z^2$ of $Z^1$ and glue both copies transversally at 
opposite ramification points. This gives the desired cover. The last assertion is obvious.
\end{proof}

Now assume that $X$ be an irreducible rational curve with 2 nodes $x_1$ and $x_2$. If $\delta_i$ denotes the number of nodes of $Z$ over $x_i$
and $r_i$ their ramification index for $i=1$ and 2, we clearly have $\delta_1 + \delta_2 = \delta$ and  
$r_i \delta_i = 6.$ We assume $\delta_1 \geq \delta_2$.  

\begin{prop} \label{prop5.2}
Let $h: Z = \cup_{i=1}^s Z_i \ra X$ be an $S_3$-cover of a rational curve $X$ with $2$ ordinary double points $x_1$
and $x_2$. Then only in the following case the cover satisfies condition $(**)$:

$s=2, \; r_1=r_2=3,\; \delta_1 = \delta_2 =2$ and the normalization of $Z$ consists of $2$ curves of genus $2$, admitting an automorphism of 
order $3$ with $4$ fixed points, 
which are glued together pairwise transversally at fixed points of $\sigma$.

There is a $1$-dimensional family of such curves. If $\tilde Y$ denotes the normalization of $Y$,
the Prym variety $P$ is an extension of the Prym variety of $J\tilde Y$ by the group $(\ZZ/3\ZZ)^2$.
\end{prop}

\begin{proof}
Let $h$ be an $S_3$-covering satisfying $(**)$. 
Suppose that $s \leq 2$. Then $c_1 = c_2 = 0$ implying $n_1 = n_2 = 0$.  Hence $r_1 = r_2 = 3$ and $\delta_1 = \delta_2 =2$.
If $s=1$, then \eqref{eq5.1} gives $g_1 = 3$. So the normalization of $Z$ is a curve of genus 3 with an $S_3$-action such that $\sigma$
admits 8 fixed points. This contradicts Proposition \ref{autgenus3}.  
If $s=2$, then \eqref{eq5.1} gives $g_1 = 2$ and we are in the case of the proposition.

If $s=3$, then $c_1 = 0, c_2 = 1$ and hence $n_1 = 0, n_2 = 1$. So, up to labeling,
we have $r_1 = 2, \delta_1 = 3$ and $r_2 = 3, \delta_2 = 2$. 
On the other hand, \eqref{eq5.1} says $9 = 3g_1 + \delta_1 + \delta_2$ which gives a contradiction.

Finally, if $s=6$, then each of the components of the normalization of $Z$ is isomorphic to $\PP^1$ and the ramification indices of 
the nodes are $r_1 = r_2 = 1$. So $c_1 = 1, \; n_1 = 2$ and $c_2=n_2 = 0$, which contradicts $(**)$. 

Concerning the existence statement, there is a one-dimensional family of smooth curves of genus 2 admitting an automorphism of order 3 (see the proof 
of proposition \ref{quot}). Each of the curves is a degree-3 cover of $\PP^1$ doubly ramified in exactly 4 points. Take 2 copies of these curves and glue 
them transversally pairwise at 2 of their ramification points. The last assertion is obvious.
\end{proof}

The remaining 4 types of non-smooth stable curves $X$ of genus 2 have 2 irreducible components which we denote by $X^1$ and $X^2$.
Let $h:Z \ra X$ be an $S_3$-cover satisfying $(**)$. For $j=1,2$ we denote
$$
Z^j := h^{-1}(X^j) = \cup_{i=1}^{s_j} Z_i^j
$$
with irreducible components $Z_i^j$. The analogue of \eqref{eq5.1} in this case is
\begin{equation}   \label{eq5.2}
7 = p_a(Z) = \sum_{j=1}^2  s_j(g_1^j -1) + \delta +1.
 \end{equation}

\begin{prop} \label{prop5.3}
Let $h: Z = \cup_{j=1}^2 \cup_{i=1}^{s_j} Z^j_i \ra X= X^1 \cup X^2$ be an $S_3$-cover with elliptic curves $X^i$ and $X^1 \cap X^2 = \{x\}$. 
Then only in the following case the cover satisfies condition $(**)$:

$s_1=s_2=1,\; r=3, \; \delta=2$. 
$Z^1$ and $Z^2$ are smooth hyperelliptic curves of genus $3$ intersecting transversally in $2$ points $z_1$ and $z_2$. 
For $j = 1,2$ the map $Z^j \ra X^j$ is an $S_3$-covering ramified 
exactly at $z_1$ and $z_2$. 

There is a $2$-dimensional family of such coverings. If $P^i$ denotes the Prym variety of the covering $Y^j \ra X^j$ for $j = 1$ and $2$, the Prym variety $P$ is isomorphic to $P^1 \times P^2$ as principally polarized abelian varieties. 
\end{prop}

\begin{proof}
Let $g_j$ denote the geometric genus of $Z^j_i$ for $j=1$ and 2.
By Proposition \ref{6components} we may assume 
$$
3 \geq s_1 \geq s_2 \geq 1.
$$

Suppose $s_1 \leq 2$. Then $c_1 = c_2 = 0$ and hence $n_1 = n_2 = 0$.  This implies $r= 3$ and thus $\delta = 2$.
If $s_1 = s_2 = 1$, Hurwitz formula implies $g_j = 3$. By Propositions \ref{autgenus3} and \ref{authypgen3}, $Z_j$ is hyperelliptic for $j = 1,2$. This gives the case of the proposition. 

Suppose $s_1 =2$. The case $s_2 =2$ cannot exist, since there is no connected graph with 4 vertices and 2 edges. 
If $s_2 = 1$, then each component $Z^1_i$ is smooth and maps 3:1 to $X^1$ with exactly one doubly ramified point. Hence it is of genus 2, which is a 
contradiction, since by Proposition \ref{quot} there is no curve of genus 2 with an automorphism of order 3 with quotient an elliptic curve. 

If $s_1 = 3, s_2 = 1$ or 2, then $c_1 = 0, c_2 = 1$ implying $n_1 = 0, n_2 =1$. This gives $r= 2$ and $\delta = 3$. But then $Z^1_j \ra X^1$ would be a 2:1 cover 
ramified exactly in one point, a contradiction.
Finally suppose $s_1 = s_2 = 3$. So $c_1=0, c_2=2$ and hence either $n_1 = 0,  n_2 = 2$ or $n_1 = 1, n_2 =0$. Both cases cannot occur, since $X$
has only one node.    

As for the existence statement,  
according to Proposition \ref{authypgen3}, there is a one-dimensional family of hyperelliptic curves $Z$ of genus 
$3$ with $S_3$-action and quotient $Z/S_3$ an elliptic curve.
Take 2 of them and intersect them transversally at the two ramification points of 
$\sigma$. The last assertion is obvious.
\end{proof}

\begin{prop} \label{prop5.4}
Let $h: Z = \cup_{j=1}^2 \cup_{i=1}^{s_j} Z^j_i \ra X= X^1 \cup X^2$ be an $S_3$-cover with a nodal rational curve $X^1$ and an elliptic 
curve $X^2$ such that
$X^1 \cap X^2 = \{x\}$. 
Then only in the following case the cover satisfies condition $(**)$:

$s_1 =2, s_2 = 1, r = r_1 = 3, \delta = \delta_1 = 2$. 
The components $Z^1_1$ and $Z^1_2$ are copies of the elliptic curve admitting an 
automorphism $\sigma$ of order $3$ with $3$ fixed points, $2$ of which are glued together transversally at opposite fixed points of $\sigma$. The third fixed point
in $Z^1_i$ for $i=1$ and $2$ is glued transversally to $Z^2$. The curve $Z^2$ is hyperelliptic of genus $3$ with an $S_3$-action, doubly 
ramified in the $2$ intersection points with $Z^1$.  

There is a $1$-dimensional family of such coverings. If $\tilde Y^1$ denotes the normalization of $Y^1$ and $P^2$ the Prym variety of the covering $Y^2 \ra Z^2$,
the Prym variety $P$ is an extension of $J \tilde Y^1 \times P^2$ by the group $\ZZ/3 \ZZ$.
\end{prop}

\begin{proof}
Let $x_1$ denote the node of $X^1$, $\delta$ (respectively $\delta_1$) the number of nodes over $x$ (respectively $x_1$) and  $r$ (respectively $r_1$) 
their ramification index. As above, let $g_j$ again denote the geometric genus of $Z^j_i$ for $j = 1,2$.

According to Proposition \ref{6components},  we may assume that $s_j =1,2$ or $3$ for $j=1,2$. In any case, the components have non-trivial stabilizer, i.e. $c_1=0$. 

If $s_j \leq 2$ for $j = 1,2$, then $c_2 = 0$ implying $n_1=n_2=0$. This gives $r=r_1=3$ and $\delta = \delta_1 =2$.
If $s_1=s_2=1$, then applying the Hurwitz formula gives $g_1=1$ and $g_2=3$. 
But there is no elliptic curve with a non-trivial $S_3$-action.

If $s_1 = 1, s_2 = 2$, then \eqref{eq5.2} implies $g_1 + 2g_2 = 5$. By the Hurwitz formula we get $g_1 = 1$, a contradiction as in the previous case.
 
If $s_1 = 2, s_2 = 1$, then \eqref{eq5.2} gives $2g_1 + g_2 = 5$.  Since
the map $Z^2 \ra X^2$ is ramified, we have $g_2 \geq 2$.  Moreover, by Proposition \ref{quot} $g_2 > 2$ and hence $g_2 = 3$. On the other hand, by the Hurwitz formula both components of $Z^1$ are of genus 1 and we are in the case of the proposition.

If $s_1 = 2,\; s_2 =3$, or $s_2=3,\; s_1=2$,  the 2 components must intersect the 3 components equally often which is impossible.

In the case $s_1=s_2=2$ we have $c_2=0$ and $n_1=n_2=0$, which implies $r_1=r=3$. Hence the 3:1 map $Z^2_i \ra X_2$
is doubly ramified in 1 point and the Hurwitz formula yields $g_2= 2$, but by Proposition \ref{quot} there is no such 
covering.

Finally, suppose $s_1=s_2=3$ then $c_2=2$, which implies that either $n_1=1, n_2=0$ or $n_1 = 0, n_2 = 2$.
The first case contradicts the fact that $Z^1$ and $Z^2$ consist
of 3 components. If $n_1=0, n_2 =1$, then $\delta = \delta_1 = 3$. So for example $Z^2_1 \ra X^2$ would be a double cover ramified at 1 point only, a contradiction. 

As for the existence statement,  
we only note that, according to Proposition \ref{authypgen3}, there is a one-dimensional family hyperelliptic curves of genus 3 with $S_3$-action and 
only finitely many possibilities for the curve $Z^1$. 
So, gluing $Z_1$ and $Z_2$ as in the proposition gives a 1-dimensional family of $S_3$-covers satisfying $(**)$. The last assertion is obvious.  
\end{proof}

\begin{prop}  \label{prop5.5}
Let $h: Z = \cup_{j=1}^2 \cup_{i=1}^{s_j} Z^j_i \ra X= X^1 \cup X^2$ be an $S_3$-cover with $X^j$ rational with $1$ node for $j=1,2$ and
$X^1 \cap X^2 = \{x\}$. 
Then only in the following case the cover satisfies condition $(**)$:

$s_1=s_2 =2$ and all ramification indices are $3$. For $i,j = 1,2$, the normalization of $Z_i^j$ is the unique elliptic curve with an automorphism of order $3$. 
The normalization of $Z^j_i \ra X^j$ has $3$ ramification points. $Z^j_i$ intersects $Z^j_{3-i}$ in $2$ of them and $Z^{3-j}_i$ in the last one.

There are only finitely many coverings of this type. If $\tilde Y^j$ denotes the normalization of $Y^j$ for $j = 1$ and $2$, the Prym variety $P$ is an extension of 
the abelian surface $\tilde Y^1 \times \tilde Y^2$ by the group $(\ZZ/3 \ZZ)^2$.
 \end{prop}                   

\begin{proof}
For $i=1,2$, let $x_i$ denote the node of $X^i$, $\delta_i$ (respectively $\delta$) the number of nodes over $x_i$ (respectively $x$) and $r_i$ (respectively $r$)
their ramification index. Let $g_j$ again denote the geometric genus of $Z^j_i$ for $j = 1,2$.

According to Proposition \ref{6components} we may assume $3 \geq s_1 \geq s_2 \geq 1$. In any case, the components have non-trivial stabilizer, i.e. $c_1=0$. 

If $s_1 \leq 2$ and $s_2 = 1$, then $c_2 =0$ and hence $n_1=n_2 =0$. So all ramification indices are 3. In particular $\delta = \delta_2 = 2$. 
Hence $Z^2 \ra \PP^1$ is a 6-fold cover doubly ramified at 6 points which gives $g_2 = 1$. But there is no elliptic curve with non-trivial
$S_3$-action.     

If $s_1 = 3 \geq s_2 \geq 1$, then in any case $\delta = 2$ or $\delta_1 = 2$. This is absurde, since $Z_1$ has 3 components and over both nodes
$x$ and $x_1$ there are at least 3 nodes.   

We are left with the case $s_1=s_2=2$.  In this situation $c_2=n_1=n_2=0$ and $r_1=r_2=r=3$. Hence $Z^j_i \ra X^j$
is a 3:1 map doubly ramified at 3 points, which gives $g_j =1$ for $j=1,2$. So the normalization of $Z^j_i$ is the unique elliptic curve with an automorphism of order 3. 
The curves $Z^1$ and $Z^2$ must be connected, since otherwise $Z$ would not be connected. This implies that all $Z^j_i$ are smooth and
the components $Z^j_1$ and $Z^j_2$ intersect transversally
in 2 points (so $p_a(Z^j)=3$) and $Z^1$, $Z^2$ intersect transversally in 2 points as well.
                                                                                
The existence as well as the last statement are obvious. 
\end{proof}

\begin{prop} \label{prop5.6}
Suppose $X= X^1 \cup X^2$ with $X^j \simeq \PP^1$ for $j=1$ and $2$ such that
$X^1 \cap X^2 = \{x_1, x_2, x_3\}$. 
There is no $S_3$-cover $h:Y \ra X$ satisfying condition $(**)$:
\end{prop}

 \begin{proof}
Suppose that $h: Z = \cup_{j=1}^2 \cup_{i=1}^{s_j} Z^j_i \ra X$ is an $S_3$-cover satisfying $(**)$ with $X$ as in the proposition. 
For $i = 1,2,3$, let $\delta_i$ denote the number of nodes over $x_i$ and $r_i$ their ramification index. 
Let $g_j$ again denote the geometric genus of $Z^j_i$ for $j = 1,2$.

According to Proposition \ref{6components} we may assume $3 \geq s_1 \geq s_2$. By Lemma \ref{smpoint}, $Z^j$ is a disjoint union of its components $Z^j_i$ for $j=1,2$, which are moreover smooth, 
since otherwise $X^j$ would have a node. 

Suppose first $2 \geq s_1 \geq s_2 = 1$. Then $c_1 = c_2 = 0$ and thus $n_1 = n_2 =0$. So $\delta_i = 2$ for all $i$ and the Hurwitz formula gives
$g_2 = 1$, contradicting the fact that an elliptic curve does not admit a non-trivial $S_3$-action.

If $s_1 =3 \geq s_2 \geq 1$, then
every component $Z^1_i$ intersects $Z^2$ in exactly 3 points $z^i_j$, where $f(z^i_j) = x_j$ and the double
cover $Z^1_j \ra X^1 \simeq \PP^1$ is ramified exactly in these 3 points. However, there is no connected 
double cover of $\PP^1$ ramified in exactly 3 points, a contradiction.

We are left with the case $s_1=s_2=2$.  
Since $\sigma$ has to map any of the components $Z^j_i$ into itself, the corresponding quotient map 
$Z^i_j \ra X^j \simeq \PP^1$
is doubly ramified in 3 points. So the Hurwitz formula implies that $Z^j_i$ is an elliptic curve for $i,j = 1,2$.
So for $j=1$ and 2 we obtain an $S_3$-cover $Z^j_1 \sqcup Z^j_2 \ra X^j$ with disjoint elliptic curves $Z^j_1$ and $Z^j_2$.
Then $\tau$ is an isomorphism of $Z^j_1$ with $Z^j_2$. For a general point $z \in Z$ we clearly have $\sigma \tau(z) \neq \tau \sigma^2(z)$.
So $Z$ is not an $S_3$-cover. This completes the proof of the proposition.
\end{proof}
 
Note that all coverings of Propositions \ref{prop5.1}, $\cdots,$ \ref{prop5.6} satisfy condition $(*)$. So we get as an immediate consequence the 
first assertion of the following Corollary.

\begin{cor}  \label{cor5.7}
Any covering $h: Z \ra X$ satisfying condition $(**)$ satisfies condition $(*)$. In particular, the Prym variety $P$ of $f:Y \ra X$ 
is a principally polarized abelian variety of dimension $2$. For any stable curve $X$ there are only finitely many covers $h$ satisfying condition
$(**)$.
\end{cor}

\begin{proof}
The last assertion follows from the proof of the propositions.
\end{proof}

 The following pictures are the dual graphs of the curves $Z$ in Propositions 5.2, 5.4 and 5.5 from left to right.

$$
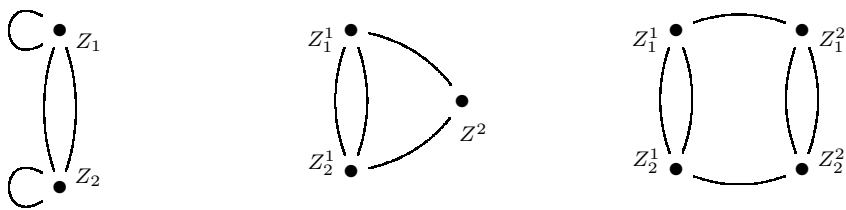
\xymatrix@R=0.7cm@C=1cm{
               {\bullet} \ar@/^/ @{-}[dd]^<{Z_1}^>{Z_2}  \ar@/_/  @{-}[dd] \ar@{-} @(ul,dl)[]   & \\
               & \\
                {\bullet} \ar@{-} @(ul,dl)[] &
}
\qquad \qquad
\xymatrix@R=.5cm@C=1cm{
               {\bullet}  \ar@/^/ @{-}[dd] \ar@/_/ @{-}[dd]_<{Z^1_1}_>{Z^1_2}  \ar@/^/ @{-}[dr] & \\
               & {\bullet} \ar@/^/  @{-}[dl]^<{Z^2}   \\
                {\bullet}  &
}
\qquad \qquad
\xymatrix@R=1.4cm@C=1.2cm{
               {\bullet}  \ar@/^/ @{-}[d] \ar@/_/ @{-}[d]_<{Z^1_1}_>{Z^1_2} \ar@/^/  @{-}[r] & 
               {\bullet}  \ar@/^/ @{-}[d]^>{Z^2_2}^<{Z^2_1}  \ar@/_/ @{-}[d] \\
                {\bullet} \ar@/_/  @{-}[r]   & {\bullet}   
}
$$

\section{Properness of the Prym map}

Let ${}_{S_3}\cM_2$ denote the moduli space of \'etale Galois covers of smooth curves of genus 2 with Galois group $S_3$ 
\cite[Theorem 17.2.11]{acg}. According to \cite{lo}, ${}_{S_3}\cM_2$ is an irreducible algebraic variety of dimension 3. A (closed) point
in it corresponds to a smooth curve $Z$ of genus 7 with an $S_3$-action and quotient $X = Z/S_3$ of genus 2. 
The quotient $Y = Z/\langle \tau \rangle$ is of genus 4, 
non-cyclic of degree 3 over $X$ and according to \cite{lo}  the Prym variety $P=P(Y/X)$ is a principally polarized abelian surface. This gives a map
$$
Pr: {}_{S_3}\cM_2 \ra \cA_2
$$
into the moduli space of principally polarized abelian surfaces, which we call the {\it Prym map}. 
In this section we show that one can extend the map to a proper map onto $\cA_2$.

Let ${}_{S_3}\overline{\cM}_2$ denote the compactification of ${}_{S_3}\cM_2$ by admissible $S_3$-covers of stable curves of genus 2 according to \cite[Chapter 17]{acg}.
To be more precise, we denote by ${}_{S_3}\overline \cM_2$ only the irreducible component containing ${}_{S_3}\cM_2$ of the moduli space as defined in \cite{acg}.
Finally, let ${}_{S_3}\widetilde{\cM}_2 \subset {}_{S_3}\overline \cM_2$ be the subset of points corresponding to covers satisfying condition $(**)$. 
The main result of this section is the following theorem.

\begin{thm} \label{thm6.1}
The set ${}_{S_3}\widetilde{\cM}_2$ is open in ${}_{S_3}\overline \cM_2$ and 
the Prym map $Pr$ extends to a proper morphism
$$
Pr: {}_{S_3}\widetilde{\cM}_2 \ra \cA_2,
$$
which is modular and which we also call the Prym map. 
\end{thm}

We mimic the proofs of \cite[6.1]{b}, \cite[I,1]{ds} and \cite[\S1]{f} using the results of 
\cite[Chapter 17]{acg}. Instead of level-$n$ structures one could also use the theory of algebraic stacks to prove the theorem.

\begin{proof}
Fix an integer $n \geq 3$ and let ${}_{S_3}\cM^{(n)}_2$ denote the moduli space of \'etale Galois covers of smooth curves 
of genus 2 with level-$n$ structure and Galois group $S_3$. As ${}_{S_3}\cM_2$, also ${}_{S_3}\cM^{(n)}_2$ is an irreducible 3-dimensional variety.
With the notation of 
\cite[Chapter 17]{acg} we have ${}_{S_3}\cM^{(n)}_2 = M_2[\psi]$, where $\psi: \pi_1(\Sigma) \ra G$ is a suitable level structure with $\Sigma$   
a fixed curve of genus 2 and $G$ a subgroup of $H^1(\Sigma, \ZZ/n\ZZ)$ with quotient $S_3$.
Let ${}_{S_3}\overline \cM^{(n)}_2$ denote the (irreducible) compactification by stable curves which is constructed in
\cite[Ch. 17, Theorem 4.8]{acg}. According to this theorem, there exists a universal family 
$\cX \ra {}_{S_3}\overline{\cM}^{(n)}_2$ of genus-2 curves with 
the corresponding structure and we have ${}_{S_3}\overline{\cM}^{(n)}_2/ Sp(4,\ZZ/n\ZZ) \simeq {}_{S_3}\overline{\cM}_2$. 
Let 
$$
\mathfrak h: \cZ \ra \cX
$$ 
denote the corresponding family of Galois covers. It is a family of admissible $S_3$-covers of stable curves of genus 7 mapping to stable 
curves of genus 2. Define the family of stable curves $\cY \ra {}_{S_3}\overline{\cM}^{(n)}_2$ by
$$
\cY := \cZ/<\tau>.
$$
The covering $\ch$ factorizes via a family non-cyclic 3-fold covers of stable curves
$$
\mathfrak f : \cY \ra \cX.
$$
According to \cite{blr}, the relative Jacobians $\Pic^0(\cY/{}_{S_3}\overline{\cM}^{(n)}_2)$ and $\Pic^0(\cX/{}_{S_3}\overline{\cM}^{(n)}_2)$ are families of semi-abelian varieties 
and the norm defines a morphism
$$
\Nm: \Pic^0(\cY/{}_{S_3}\overline{\cM}^{(n)}_2)   \ra \Pic^0(\cX/{}_{S_3}\overline{\cM}^{(n)}_2).
$$
Define the family of Prym varieties $\cP \ra {}_{S_3}\cM^{(n)}_2$ of $\mathfrak f$ by
$$
 \cP := \Ker (\Nm: \Pic^0(\cY/{}_{S_3}\overline{\cM}^{(n)}_2)   \ra \Pic^0(\cX/{}_{S_3}\overline{\cM}^{(n)}_2)).
$$                                                                                                                                   
Let ${}_{S_3}\widetilde \cM^{(n)}_2$ denote the set consisting of points  $ s \in {}_{S_3}\overline{\cM}^{(n)}_2$ such that the fibre $\cP_s$ is an abelian variety. 
Hence ${}_{S_3}\widetilde \cM^{(n)}_2$ coincides with the 
subset over which the map $\cP \ra {}_{S_3}\overline{\cM}^{(n)}_2$ is proper. This implies also that ${}_{S_3}\widetilde \cM^{(n)}_2$ 
is an open set in ${}_{S_3}\overline{\cM}^{(n)}_2$. By abuse of notation we denote by $\cP$ also the 
restriction of $\cP$ to ${}_{S_3}\widetilde \cM^{(n)}_2$. According to Proposition \ref{lem3.1} the set ${}_{S_3}\widetilde \cM^{(n)}_2$ 
coincides with the set of points $s$ such that the cover $\ch_s: \cZ_s \ra \cX_s$ satisfies condition $(**)$ and hence, by Corollary \ref{cor5.7}, 
satisfies condition $(*)$. By Corollary \ref{cor2.5}, for every $s \in {}_{S_3}\widetilde \cM^{(n)}_2$ the abelian variety $\cP_s$ admits a principal polarization $\Xi_s$.
These principal polarizations glue together in the usual way to give a family of principal polarizations, i.e. an isomorphism 
$ \phi_{\Xi}: \cP \ra \widehat{\cP}$ over ${}_{S_3}\widetilde{\cM}^{(n)}_2$. Since all families occuring are flat, we get a flat family of principally polarized abelian 
varieties of dimension 2 over ${}_{S_3}\widetilde \cM^{(n)}_2$. Hence we get a morphism
$$
p: {}_{S_3}\widetilde \cM^{(n)}_2 \ra \cA_2
$$
into the moduli space $\cA_2$ of principally polarized abelian surfaces. We follow the lines of the proof of \cite[Proposition 6.3]{b} 
to show that the map $p$ is proper. By means of the valuative criterion of properness and the completeness of ${}_{S_3}\overline{\cM}^{(n)}_2$
it is enough to show the following:
Let $T$ be the spectrum of a complete discrete valuation ring and $\eta$ its generic point. 
Consider $\tilde{\cZ} \ra T$ a family of  admissible $S_3$-coverings such that the covering $\tilde{\cZ_{\eta}} \ra \tilde{\cZ_{\eta}}/S_3$
satisfies the condition $(**)$ (and then satisfies $(*)$) and the corresponding Prym variety $\cP_{\eta}$ extends to an abelian variety over $T$. Then $\cP_s$ is an abelian variety. 
Since $\cP_s$ is an extension of an abelian variety by a torus it is isomorphic to the neutral component of the Neron model of 
$\cP_{\eta}$ over $T$ which is abelian by hypothesis. This proves the properness of $p$.

Now ${}_{S_3}\widetilde \cM^{(n)}_2 \subset {}_{S_3}\overline{\cM}^{(n)}_2$ is stable under the action of the group $Sp(4,\ZZ/n\ZZ)$. Hence the quotient
$$
{}_{S_3}\widetilde \cM_2 := {}_{S_3}\widetilde \cM^{(n)}_2/Sp(4,\ZZ/n\ZZ)
$$
is a non-empty open set in ${}_{S_3}\overline{\cM}_2$. Moreover, since $p$ commutes with this action, we obtain an induced map
$$
Pr: {}_{S_3}\widetilde \cM_2 \ra \cA_2.
$$
This is the (extended) Prym map of the theorem. The induced map $Pr$ is proper, since $p$ is. Finally, its restriction to the open set of smooth covers in 
${}_{S_3}\widetilde \cM_2$
clearly is the Prym map of \cite{lo}. 
\end{proof}

As an immediate consequence of Theorem \ref{thm6.1} we obtain

 \begin{cor} \label{cor6.2}
The extended Prym map $Pr: {}_{S_3}\widetilde{\cM}_2 \ra \cA_2$ is surjective. In other words,
every principally polarized abelian surface occurs as the Prym variety of a non-cyclic degree-$3$ admissible cover $f: Y \ra X$ of a stable curve $X$ of genus $2$.   
 \end{cor}

\section{Stratification of ${}_{S_3}\widetilde{\cM}_2$}

Consider the following stratification of the moduli space ${}_{S_3}\widetilde{\cM}_2$:
$$
{}_{S_3}\widetilde{\cM}_2 = {}_{S_3}{\cM}_2 \;\; \sqcup \;\; R \;\; \sqcup \;\; S
$$
with boundary components $R$ and $S$ where
$$ 
R = \cup_{i=1}^2 R_i \quad \mbox{and} \quad S = \cup_{i=0}^2 S_i.
$$
Here $R_2$, respectively $S_2$, denotes the 2-dimensional subspace of ${}_{S_3}\widetilde{\cM}_2$ of Proposition \ref{prop5.1}, respectively \ref{prop5.3}, 
$R_1$ respectively $S_1$, denotes the 1-dimensional subspace of ${}_{S_3}\widetilde{\cM}_2$ of Proposition \ref{prop5.2} respectively \ref{prop5.4} and finally
$S_0$ denotes the 0-dimensional subspace of ${}_{S_3}\widetilde{\cM}_2$ of Proposition \ref{prop5.5}.
Note that $R_{1}$ is in the closure of $R_2$ and $S_{i-1}$ is in the closure of $S_{i}$ for $i =1$ and 2. 
Note moreover that $S$ is closed in ${}_{S_3}\widetilde{\cM}_2$ whereas we have for the closure of $R$,
$$
\overline{R} = R \cup S_1 \cup S_0.
$$

In this section we will study the images of $R_i$ and $S_i$
under the extended Prym map $Pr$. First we determine the image $Pr({}_{S_3}{\cM}_2)$ of the open set of smooth covers. 
We will use the following well known fact: Every principally polarized abelian surface is either the Jacobian of a smooth curve of genus 2 or a canonically 
polarized product of 2 elliptic curves. We denote by $\cJ_2 \subset \cA_2$ the (open) subset of Jacobians of smooth curves and by $\cE_2$ its complement in
$\cA_2$.

Given a smooth genus-2 curve $\Sigma$ we denote by $\varphi: \Sigma \ra \PP^1$ the corresponding hyperelliptic cover. 
For any 3 Weierstrass points $w_1,w_2,w_3$ of $\Sigma$ let
$\varphi_{2(w_1+w_2+w_3)}$ denote the map $\Sigma \ra \PP^1$ defined by the pencil $(\lambda (2(w_1+w_2+w_3)) + \mu (2(w_4+w_5+w_6))_{(\lambda,\mu) \in \PP^1}$ where $w_4,w_5,w_6$ 
are the complementary Weierstrass points. 
According to 
\cite[Proof of Theorem 5.1]{lo}, the curve
$\Sigma$ fits into the following commutative diagram
$$
\xymatrix{
Y \ar[r] \ar[d]_{f}  & \PP^1\ar[d]^{\bar{f}} & \ar[l]_{\varphi} \Sigma \ar[dl]^{\psi} \\
X \ar[r] & \PP^1 &
    }
$$ 
where the horizontal maps are the hyperelliptic covers, and $\psi$ is given by a pencil $g^1_6 \subset |3K_{\Sigma}|$ 
which, up to enumeration, is generated by the 2 divisors $2w_1 + 2w_2 + 2w_3$ and  $2w_4 + 2w_5 + 2w_6$. 


\begin{prop} \label{prop7.1}
$Pr ({}_{S_3}\cM_2) =\{  J\Sigma \in \cJ_2 \mid  \  \exists \; w_1,w_2,w_3 \mbox{ Weierstrass points of }  \Sigma
\mbox{ such } $
 
$\hspace{5.5cm} \mbox{  that } \varphi_{2(w_1+w_2 +w_3)} = \bar{f} \circ \varphi, \mbox {and $\bar{f}$ is simply ramified } \}$
\end{prop}

\begin{proof}
The inclusion ''$\supset$'' has been shown in \cite{lo}. Conversely, let $J\Sigma$ be the Prym variety of the cover $f:Y \ra X$ 
associated to an element of $ {}_{S_3}\cM_2$. According to the diagram above, for a choice of 3 Weierstrass points 
$w_1, w_2, w_3$ of $\Sigma$ the map $\varphi_{2(w_1+ w_2+ w_3)}: \Sigma \ra \PP^1$ factorizes through  the hyperelliptic covering
$\varphi$ and a map $\bar{f}: \PP^1 \ra \PP^1$.  Suppose that $\bar{f}$ is not simply ramified, then the branch locus of $\psi$ consists of at  most 5 points. 
Two of these branch points are images of $w_1, w_2, w_3$ and $w_4, w_5, w_6$ and the other are the branch locus of $\bar{f}$.  
Hence there exists a point $q \in \PP^1$ outside of branch locus of $\psi$ such that ${f}^{-1}(q)$ contains the image 
of a Weierstrass point of $Y$. Since $\bar{f}$ is \'etale on the fiber of $q$, all the 3 points in $\bar{f}^{-1}(q)$ are images of  Weierstrass points
 of $Y$, contradicting the number of Weierstrass points of $Y$.
\end{proof} 

\begin{prop} \label{prop7.2}
The extended Prym map restricts to a surjective morphism (denoted by the same letter)
$$
Pr: \sqcup_{i=0}^2 S_i \ra \cE_2.
$$
\end{prop}

\begin{proof}
For $i = 2$ we saw that $Pr(S_2) \subset \cE_2$ already in Proposition \ref{prop5.3}. This implies the analogous statement for $i = 1$ 
and 2, since any specialization of a principally polarized product of elliptic curves is itself such a product. Moreover, 
$\sqcup_{i=0}^2 S_i$ is closed in ${}_{S_3}\widetilde{\cM}_2$. So $Pr: \sqcup_{i=0}^2 S_i \ra \cE_2$ is a proper morphism 
by Theorem \ref{thm6.1}. Since it is clearly dominant, the surjectivity follows from this.
\end{proof}

Finally, we determine the image of $R$ under the Prym map. Consider first $R_1$.

\begin{prop}  \label{prop7.3}
 $\Pr(R_1) \subset \cJ_2$.
\end{prop}
 
\begin{proof}
Let $f: Y \ra X$ be the cover given by an element of $R_1$. So $Y$ is  an irreducible curve of geometric 
genus 2 with two nodes and $X$ is a rational irreducible curve with 2 nodes. Let $\tY$ and $\tX$ be their normalizations. 
Then   $Y= \tilde Y/ (\tilde y_1 \sim \tilde y_2, \tilde y_3 \sim \tilde y_4)$, where $\tilde y_i$ for $i=1, \ldots , 4$, are doubly ramified 
points of the 
triple covering $\tilde f: \tY \ra \tX$ and $X = \tX / (\tilde x_1 \sim \tilde x_2, \tilde x_3\sim \tilde x_4)$, with $\tilde x_i = \tilde f (\tilde y_i )$ for $i = 1,\ldots,4$. 
According to Proposition \ref{prop5.2}, the associated Prym variety $P$ is an extension
$$
0 \ra (\ZZ/ 3\ZZ)^2 \ra P \stackrel{\pi}{\ra}  J\tY \ra 0.
$$

Suppose that $P \simeq E_1 \times E_2 $ as principally polarized abelian surfaces, with $E_1, E_2$ elliptic curves and $E_1 \times E_2$ 
with canonical (split) polarization. 
The pull-back $\pi^*(\Theta_{\tY}) $ of the canonical principal polarization of $J\tY$ 
defines a covering of degree 9, $\pi^*(\Theta_{\tY}) \ra \Theta_{\tY} $, so it contains an irreducible component of genus at least 2
(otherwise the map $\pi$ would have fibres of positive dimension, which is impossible). 
On the other hand, according to the assumption and the construction in Proposition \ref{prop5.2}, $\pi^*(\Theta_{\tY})$ defines 
the 3-fold of the canonical principal polarization on $E_1 \times E_2$. Hence the linear system of this polarization contains an irreducible 
curve of genus at least 2.
This contradicts the K\"unneth formula, according to which all (reduced) irreducible components of this linear system are elliptic curves.  
 \end{proof}

As an immediate consequence of Proposition \ref{prop7.3}, we get that the Prym variety of a general element if $R_2$ is the Jacobian of a genus 2 curve. We will see
in the next section that this is the case for every element in $R_2$.

According to Proposition \ref{prop7.1}, the set
$$
\cJ_2^u := Pr({}_{S_3}\cM_2) \subset \cJ_2 
$$ 
is the set of Jacobians which admit Weierstrass points $w_1,w_2,w_3$ such that the map $\varphi_{2(w_1+w_2+w_3)}$ factors via the hyperelliptic cover and a 
simply ramified map $\bar f: \PP^1 \ra \PP^1$. It is easy to see that $\cJ_2^u$ is open in $\cJ_2$ and thus in $\cA_2$. Denote by
$$
\cJ_2^r:=\{  J\Sigma \in \cJ_2 \mid  \  \exists \; w_1,w_2,w_3 \mbox{ Weierstrass points in } \Sigma
\mbox{ with } \bar{f} \mbox{ not simply ramified} \}.
$$
Note that $\cJ_2^u \cap \cJ_2^r \neq \emptyset$.
So the covering
$$
\cA_2 = \cJ_2^u \cup \cJ_2^r \sqcup \cE_2.
$$
is not a stratification of $\cA_2$.
The following theorem summarizes the situation

\begin{thm} \label{imagethm}
The stratification ${}_{S_3}\widetilde{\cM}_2 = {}_{S_3}{\cM}_2 \; \sqcup \; R \; \sqcup \; S$ and the covering
$\cA_2 = \cJ_2^u \cup \cJ_2^r \sqcup \cE_2$ are compatible under the Prym map $Pr$. To be more precise, we have
\begin{enumerate}
 \item  $Pr({}_{S_3}{\cM}_2) = \cJ_2^u,$

\item $Pr(R) \subset \cJ^r_2$ and

\item $Pr (S)  = \cE_2 $. 
\end{enumerate}
\end{thm}

\begin{proof}
Part (1) has been shown in Proposition \ref{prop7.1} and part (3) in Proposition \ref{prop7.2}. 
Concerning (2),  $Pr(R_1) \subset \cJ^r_2$ according to Proposition \ref{prop7.3} and Remark \ref{2nodes}.
In  Proposition \ref{prop8.6} we will show that $Pr(R_2) \subset \cJ_2$. Let $f: Y \ra X$ be a covering given by an element in $R_2$,
then $Y= \tY / \tilde{y_1} \sim \tilde{y_2}$  where $\tY$ is the normalization of $Y$. It is not difficult to see that, for a given 
map $\varphi_{2(w_1+w_2+w_3)}$ (see diagram (\ref{commdiag})), the corresponding $\bar{f} : \PP^1 \ra \PP^1$ is 
doubly ramified  in the image of the points $\tilde y_1$ and $\tilde y_2$. This shows that $Pr(R_2) \subset \cJ^r_2$.

\end{proof}

\section{The image of $R_2$ under the Prym map}

Suppose that we are given an $S_3$-cover in $R_2$, which gives the cover $f: Y \ra X$.
So $X$ is an irreducible curve of geometric genus 1 with one node and normalization $\widetilde X$, i.e. 
$$
X = \tilde X/\tilde x_1 \sim \tilde x_2,
$$  
with points $\tilde x_1 \neq \tilde x_2$ of $\tilde X$. The curve $Y$ is an irreducible curve of geometric genus 3 with one node 
$y_0$ and normalization $\tilde Y$, i.e. 
$$
Y = \tilde Y/\tilde y_1 \sim \tilde y_2,                                          
$$ 
with points $\tilde y_1 \neq \tilde y_2$ of $\tilde Y$ such that 
$\tilde f: \tilde Y \ra \tilde X$ is doubly ramified exactly at $\tilde y_1$ and $\tilde y_2$. As a degeneration of hyperelliptic curves, $Y$ is hyperelliptic and hence $\tilde Y$ is hyperelliptic.
As an elliptic curve, $\tilde X$ admits a one-dimensional family of double covers of $\PP^1$. There is exactly one such double 
cover such that the square in the following diagram commutes. 
\begin{equation} \label{diagr}
\xymatrix{
\tY \ar[rr]^{h_{\tilde Y}} \ar[dd]_{\tilde f} \ar[dr]_{n_{\tilde Y}} && \PP^1\ar[dd]^{\bar{f}}  \\
& Y  \ar[ur]_{h_Y} \ar[dd]_(.3){f}  | ! {[d] } \hole & \\
\tX \ar[rr]^{h_{\tX} \qquad} \ar[dr]_{n_Y} && \PP^1 \\
& X \ar[ur]_{h_X}&
    }
\end{equation}
In particular $f$  and $\tilde{f}$ commute with the respective hyperelliptic involutions. 
We denote the map $\tilde X \ra \PP^1$ by $h_{\tilde X}$ and call it (by abuse of notation) the hyperelliptic cover of the elliptic curve $\tilde X$. 
The hyperelliptic involution of $\tilde Y$ exchanges the points $\tilde y_1$ and $\tilde y_2$ and the corresponding involution 
on $\tilde X$ exchanges the points $\tilde x_1$ and $\tilde x_2$. 
Hence $\bar{f}$ is doubly ramified at $h_{\tilde Y}(\tilde y_1) = h_{\tilde Y}(\tilde y_2)$ and simply ramified at 2 other points. 

\begin{rem}\label{2nodes}
For a $S_3$-covering $f:Y\ra X$ given by an element in $R_1$, i.e.  a covering with two nodes as in Proposition \ref{prop7.3},
the corresponding maps $\tilde{f}: \tY \ra \tX $ and $\bar{f}: \PP^1 \ra \PP^1$ fit into a commutative diagram  as in \ref{diagr}
and it shows that $\bar{f}$ is doubly ramified at the points $h_{\tilde Y}(\tilde y_1) = h_{\tilde Y}(\tilde y_2)$
and $h_{\tilde Y}(\tilde y_3) = h_{\tilde Y}(\tilde y_4)$. According to Proposition \ref{prop7.3},  the Prym variety of $f$ is 
isomorphic to the Jacobian of a genus 2 curve $\Sigma$, which possesses $w_1, w_2, w_3$ Weierstrass points such that 
$\phi_{2(w_1+w_2+w_3)} = \bar{f} \circ \phi$. Then $Pr(f) \in \cJ_2^r$.
\end{rem}

According to \cite[Section 2]{b}, there is a canonical theta divisor in $\Pic^3 (Y)$, namely
$$
\Theta_Y := \{ L \in \Pic^3(Y) \mid  h^0(Y, L) \geq 1 \}.
$$
The following proposition describes $\Theta_Y$ set-theoretically. For this we need some notation.
As $n_Y: \tY \ra Y$  denotes the normalization of $Y$, we have a surjective morphism
$$
\Pic^3 (Y) \stackrel{n^*_Y}{\ra} \Pic^3 (\tilde Y) .
$$
For any $M \in \Pic^3 (\tilde Y)$ we denote by 
$$
F(M) := (n^*_Y)^{-1}(M) \simeq \CC^*,
$$
the set of line bundles $L$ in $\Pic^3(Y)$ mapping to $M$. If we fix an isomorphism 
$$
\varphi: M_{\tilde y_1} \stackrel{\simeq}{\ra} M_{\tilde y_2},
$$ 
then $L$ is determined by a constant $c \in \CC^*$
such that $M_{\tilde y_1}$ is glued with $M_{\tilde y_2}$ by $c \varphi$. We denote by
$$
\Theta(M) := F(M) \cap \Theta_Y,
$$
the fibre over $M$ of the restricted map $n^*_Y|_{\Theta_Y}: \Theta_Y \ra \Pic^3 (\tilde Y)$. 
In the sequel we denote 
$$Y^0 := Y \setminus \{y_0\},
$$
the smooth locus of $Y$. 
The {\it Abel map} in degree 3 of $Y$ is defined as the morphism
$$
\alpha^3_Y: (Y^0)^3 \ra \Pic^3 (Y), \qquad (y_1,y_2,y_3) \mapsto \cO_Y(y_1 + y_2 + y_3).
$$
Let $H_{\tY}$ denote the hyperelliptic line bundle of $\tilde Y$. Since the line bundles of degree 3 on $\tilde Y$ with $h^0 = 2$ are exactly the line bundles
$H_{\tilde Y}(y)$ for some point $y \in \tilde Y$ and $y$ is a base point of the corresponding linear system, we have for $M \in \Pic^3(\tilde Y)$,
$$
h^0(M) = 2  \quad \Leftrightarrow \quad M \simeq H_{\tilde Y}(y) \; \mbox{for some}\; y \in \tilde Y,
$$
$$
h^0(M) = 1 \quad \Leftrightarrow \quad  M \simeq \cO_{\tilde Y}(y_1 + y_2 + y_3) \; \mbox{with} \; M(-y_i) \not \simeq H_{\tilde Y} \; \mbox{for} \; i=1,2,3.
$$
The following proposition is easy to prove, but since it is a special case of 2 more general lemmas of \cite{c}, we refer to this paper instead.

\begin{prop} \label{prop8.1}
For any $M \in \Pic^3 (\tilde Y)$ exactly one of the following cases occurs,
\begin{enumerate}
\item If $M \simeq H_{\tilde Y}(y)$ for some $y \in \tilde Y$, $y \neq \tilde y_1,\tilde y_2$, then
$$
\Theta(M) = F(M) \; \mbox{and there is a unique} \; L_M \in \Theta(M) \; \mbox{with} \; L_M \in \im(\alpha^3_Y).
$$ 
Moreover,
$h^0(L_M) =2$ and $h^0(L) = 1$ for all $L \in \Theta(M), \; L \not \simeq L_M$. 
\item If $M \simeq H_{\tilde Y}(\tilde y_i)$ for $i = 1$ or $2$, then
$$
\Theta(M) = F(M) 
$$
and $ F(M) \cap \im(\alpha^3_Y) = \emptyset.$
Moreover, $h^0(L) =1$ for all $L \in \Theta(M)$.
\item If $M \simeq \cO_{\tilde Y}(y_1 + y_2 + y_3)$ with $M(-y_i) \not \simeq H_{\tilde Y}$,  then $\Theta(M)$ contains exactly one line bundle denoted $L_M$ and
 $$
\Theta(M) = \{L_M\} 
$$
with $L_M \in \im (\alpha^3_Y).$ Moreover, $h^0(L_M) = 1$.
\item If $M \simeq \cO(y_1 + y_2 + \tilde y_i)$ for $i = 1$ or $2$ with $M(-\tilde y_i) \not \simeq H_{\tilde Y}$ and $y_j \neq \tilde y_{2-j}$ for $j=1$ and $2$, then
$$
\Theta(M) = \emptyset.
$$
\end{enumerate}
 \end{prop}

\begin{proof}
We have $h^0(M) \geq 1$, since $\tilde Y$ is of genus 3, and $h^0(M) \leq 2$ by Clifford's theorem. Hence exactly one of the 4 cases occurs.
Assertion (1) is a special case of \cite[Lemma 2.2.4(2)]{c}, (2) is a special case of \cite[Lemma 2.2.4(3)]{c}, (3) a special cases of \cite[Lemma 2.2.3(2a)]{c}
and finally (4) a special case of \cite[Lemma 2.2.3(1)]{c}. 
\end{proof}

\begin{rem}
If $M = \cO_{\tY}(y_1 + y_2  + y_3)$ is as in (1) or (3) in the proposition, then we can choose $y_1, y_2, y_3$ in such a way that they correspond to points of $Y^0$, 
which we denote by the same letters. Then the uniquely determined line bundle $L_M$ of Proposition \ref{prop8.1} is just $\cO_Y(y_1 + y_2 + y_3)$.
We have for any $\cO_{\tY}(y_1 + y_2 + y_3), \cO_{\tY}(y'_1 + y'_2 + y'_3)$ of type (1) or (3),
\begin{equation}
\cO_Y(y_1 + y_2 + y_3) \simeq \cO_Y(y'_1 + y'_2 + y'_3)  \Leftrightarrow y_1 + y_2 + y_3 \sim_{\tY} y'_1 + y'_2 + y'_3,
\end{equation}
where $\sim_{\tY}$ denotes linear equivalence on $\tY$. 
 \end{rem}

For the rest of this section we fix the following notation.
Let $q_1, \ldots, q_8$ denote the Weierstrass points of $\tilde Y$ and $p_1, \ldots, p_4$ the ramification points of $h_{\tilde X}$.
From the commutativity of the square in the diagram \eqref{diagr} $\tilde{f}$ maps the Weierstrass points of $\tY$ to the Weierrstrass points of $\tX$ and the hyperelliptic involution of $\tY$ acts on each fibre of $\tilde{f}$. Thus the fibre over a point   
$p_i$, $i=1, \dots, 4$ consists of  1 or 3 Weierstrass points. We conclude that there are 2 points $p_j$ such that $\tilde f^{-1}(p_j)$ consists of 3 of the points $q_i$ and for the remaining 2, $\tilde f^{-1}(p_j)$ contains
of 1 of the points $q_j$. Without loss of generality, we may assume that
$$
\tilde f(q_1)  = \tilde f(q_2) = \tilde f(q_3) = p_1, \quad \tilde f(q_4)  = \tilde f(q_5) = \tilde f(q_6) = p_2
$$
and 
$$
\tilde f(q_7) = p_3, \quad \tilde f(q_8) = p_4.
$$
If $H_Y$ and $H_X$ denote the hyperelliptic line bundles of $Y$ and $X$, i.e. the line bundles defining the hyperelliptic covers $h_Y$ and $h_X$, we have
$$
f^*(H_X) \simeq H_Y^3.
$$ 

In order to describe the restriction of the theta divisor of $JY$ to the Prym variety $P$ we will work in $\Pic^2(Y)$. For this we consider 
the following translate of $\Theta_Y$:
$$
\Theta_{q_1} := \Theta_Y - q_1 \subset \Pic^2(Y).
$$
A priori $\Pic^2(X), \Pic^2(Y), \Pic^2(\tX),$ and  $\Pic^2(\tY)$ are not algebraic groups, but they have a canonical point, namely the hyperelliptic line bundle.
Using this, we may consider them as semiabelian varieties and have the following commutative diagram with exact rows,
\begin{equation} \label{deg2diag}
\xymatrix{
0 \ar[r] & \ZZ/3\ZZ \ar[r] \ar[d] & P \ar[r] \ar[d] & \tP \ar[r] \ar[d] & 0\\
1 \ar[r] & \CC^* \ar[r]^{i_Y} \ar[d]^{Nm_{\CC^*}} & \Pic^2(Y) \ar[r]^{n_Y^*} \ar[d]^{Nm_f} & \Pic^2(\tY) \ar[r] \ar[d]^{Nm_{\tilde f}} & 0\\
1  \ar[r] & \CC^* \ar[r]^{i_X}  & \Pic^2(X) \ar[r]^{n_X^*}  & \Pic^2(\tX) \ar[r]  & 0\\
}
\end{equation}
where $Nm_{\CC^{*}}$ is the third power map. The maps $i_Y$ and $i_X$ are defined as follows: Recall that $H_Y$ is defined in terms of $H_{\tY}$ by
gluing the fibres at $\tilde y_1$ and $\tilde y_2$. The constant $c \in \CC^*$ corresponding to $H_Y$ depends on the gluing. We glue the fibres of $H_{\tY}$ in such a way that
the constant corresponding to $H_Y$ is 1. In other words we have $H_Y = i_Y(1)$ and similarly for $H_X$. So $P= \Nm_f^{-1}(H_X)$ and 
$\tP= \Nm_{\tilde f}^{-1}(H_{\tX})$.

For the other fibres we proceed similarly: Given $M \in \Pic^2(\tY)$, we choose a line bundle $L_M \in n_Y^*(M)$ and glue the fibres $M_{\tilde y_1}$ and 
$M_{\tilde y_2}$ in such a way 
that the constant corresponding to $L_M$ is 1. Then we have
$$
F(M) = \{ L_M \otimes i_Y(c) \otimes H_Y^{-1} \;|\; c \in \CC^* \}.
$$   
If $M$ is of type (1) or (3) of Proposition \ref{prop8.1}, we choose $L_M$ as described there.

\begin{prop} \label{Xi}
A line bundle $L\in \Pic^2(Y)$ is in the intersection $P \cap \Theta_{q_1}$ if and only if 

{\em (a)} $L \simeq \cO_{Y} (y_1 +y_2 +q_i- q_1) )$ with $f(\iota_Y y_2 ) = f(y_2)$  for $i=1,2,3$ or 

{\em (b)} $L\simeq \cO_Y(q_i -  q_1) \otimes i_Y(\zeta^j) $ where $ \zeta$ is a primitive third root of the unity and $i = 1,2$ or $3$,  
$j=0,1$ or $2$.

\end{prop}

\begin{proof}
The elements in $\Theta_Y$ have been described in Proposition \ref{prop8.1}.  Let $L \in \Theta_{q_1}$ and set $M= n_Y^*(L(q_1))$. 
Suppose $M$ is of type  $(1)$ of Proposition \ref{prop8.1}, i.e. $M\simeq H_{\tY}(y)$ for some $y \in \tY \setminus \{\tilde y_1, \tilde y_2\}$. Then 
$$
L = \cO_Y(y- q_1) \otimes i_Y(c)  \quad \mbox{for some} \quad c \in \CC^*. 
$$
So $\Nm_f(L) = \cO_X( f(y)- p_1 ) \otimes i_X(c^3)$.
Since an element of $\Pic^2(X)$ is uniquely determined by its image in $\Pic^2(\tX)$ and an element in $\im i_X$,   
$\Nm_f(L)=H_X$ if and only if $i_X(c^3) = H_X$ and $f(y)=p_1$, that is, $c^3=1$ and
$y \in f^{-1}(p_1)=\{q_1, q_2, q_3\}$ .  This gives the line bundles in (b).

Suppose $M\simeq H_{\tY}(\tilde y_i)$ for $i=1$ or 2.  If $\Nm_f(L)= H_X$, by the commutativity of diagram \eqref{deg2diag}, 
$\Nm_{\tilde f}(M(-q_1)) = H_{\tX}(\tilde x_i - p_1) = H_{\tX}$, which implies $\tilde x_i=p_1$, a contradiction.

Suppose now that $M$ is of type (3), i.e. $M\simeq \cO_{\tY} (y_1+y_2+y_3)$ with $M(-y_i) \not \simeq H_{\tY}$ for all $i$. Then 
$\Theta(M) $ consists only of the line bundle $L_M\in \im (\alpha^3_Y)$. So 
$$
\Nm_f(L_M(-q_1))=   \Nm_f(\cO(y_1 +y_2  +y_3 -q_1)) = \cO_X(f(y_1)+f(y_2) + f(y_3)- p_1).
$$
and  $\Nm_f(L_M(-q_1))= H_X$ if and only if 
$$
\cO_X(f(y_1)+f(y_2) + f(y_3)) \simeq  H_X \otimes \cO_Y(p_1).
$$ 
Since $p_1$ is a base point of the linear system $|H_X \otimes \cO_Y(p_1)|$ we
conclude that, after possibly renumerating, we get $y_3\in f^{-1}(p_1)$ and $f(y_1)+f(y_2) \sim H_X  $. This gives the line bundles in (a).
\end{proof}

For $i=1,2,3$, we consider the following sets 
\begin{eqnarray*}
\widetilde{\Xi}_i & := & \{ \cO_{Y}( y+ z +q_i -q_1) \in \Pic^2(Y) \mid  y, z \neq \tilde y_1,\tilde y_2, \;  f(\iota_Y z) =f(y) \} \\
& &  \cup \quad  \{ \cO_Y(q_i -  q_1) \otimes i_Y(\zeta^j) \mid j=0,1,2 \},
\end{eqnarray*}
As a consequence of Proposition \ref{Xi} we obtain, 
$$
P \cap \Theta_{q_1} =\widetilde{\Xi}_1 \cup \widetilde{\Xi}_2 \cup \widetilde{\Xi}_3.
$$
For $i =1,2,3$, define the scheme $ \Xi_i$ as the closure of $ \widetilde{\Xi}_i \setminus  \{ \cO_Y(q_i -  q_1) \otimes i_Y(\zeta^j) \mid j=0,1,2 \} $
with reduced subscheme structure.
\begin{lem}
The scheme $\Xi_i $ is a complete curve for $i=1,2,3$ and we have the following equality of sets
$$
P \cap \Theta_{q_1} = \Xi_1 \cup \Xi_2 \cup \Xi_3.
$$
\end{lem}
\begin{proof}
Note first that $\Xi_i = \Xi_1 +q_i -q_1$ for $i = 1$ and 2. So for the proof we have only to show the assertion for $i =1$ and for this it suffices to show that
$P \not \in \Theta_{q_1}$, since then $P \cap \Theta_{q_1}$ is a divisor in a surface. 
The proof of this is very similar to the smooth case as given in \cite[Lemma 4.11]{lo} and we omit it. 

For the last assertion, certainly the right hand side is contained in the left hand side and the difference consists of at most the finitely many points
$\{ \cO_Y(q_i -  q_1) \otimes i_Y(\zeta^j) \mid j=0,1,2 \}$. But $P \cap \Theta_{q_1}$ is a divisor in $P$, so does not contain isolated points. 
So these points are contained in the right hand side.
Together with 
the first assertion this completes the proof of the lemma.
\end{proof}

\begin{prop} \label{prop8.5}
The principal polarization $\Xi$ of the Prym variety $P$ is given by each of the algebraically equivalent divisors $\Xi_1 \equiv
\Xi_2 \equiv \Xi_3$. 
\end{prop}

\begin{proof}
The proof is the same as for \cite[Theorem 4.13]{lo} and will be omitted.
\end{proof}

Now we are in a position to prove the main result of this section, which completes the proof of Theorem \ref{imagethm}.

\begin{prop} \label{prop8.6}
The Prym variety of any element of $R_2$ is the Jacobian of a  smooth curve of genus $2$.
\end{prop}

\begin{proof} 
Let $f: Y \ra X$ be the cover associated to an element of $R_2$ and let 
$\Xi := \Xi_1$ denote the curve which, according to Proposition \ref{prop8.5}, defines a principal polarization of the corresponding Prym variety $P$. 
Then $\Xi$ is either a smooth genus 2 curve or the union of two elliptic curves intersecting transversally in one point.  

Suppose $\Xi$ is reducible. Then also the (non-complete) curve  
$$
\Xi^{0} :=  \Xi  \setminus  \{ i_Y(\zeta^j) \mid j=0,1,2 \} 
$$ 
is reducible.  The curve $\Xi^{0}$ is isomorphic to  its image in $\Pic^2(\tY)$, denoted also by $\Xi^{0}$.  Since 
$H_{\tY} \notin \Xi^{0}$, we can and will consider $\Xi^{0}$ as a
subset in the symmetric product $\tY^{(2)}$.  

Away from the points $\tilde y_1$ and $\tilde y_2$ the map $\tilde f: \tY \ra \tX$ is \'etale. So for every point $y \in \tY, \; y \neq \tilde y_1,\tilde y_2$ the fibre 
$\tilde f^{-1}\tilde f(y)$
consists of exactly 3 points. Let us denote these points by $\tilde f^{-1}\tilde f(y) = \{y,y',y''\}$; we denote $\tY^{0} := \tY \setminus \{\tilde y_1,\tilde y_2\}$ and 
define the curve 
$$
D:=\{ (y, \iota_{\tY} y'), (y, \iota_{\tY} y'') \in \tY \times \tY \mid y \in \tY \setminus \{\tilde y_1,\tilde y_2\} \}
$$
with reduced scheme structure and denote by $\bar D$ its completion.
The restriction of the canonical map $\tY \times \tY \ra \tY^{(2)}$ to $D$ defines a $2:1$ map $D \ra  \Xi^{0}$, which extends to the closure 
$\bar{D} \ra \bar{\Xi}^{0} = \Xi$. 

On the other hand, the projection to the first component gives a $2:1$ map  $\bar{D} \ra \tY$, which is ramified exactly at the points 
$(\tilde y_1, \iota_{\tY} \tilde y_1) = (\tilde y_1,\tilde y_2)$ and  $(\tilde y_2, \iota_{\tY} \tilde y_2) = (\tilde y_2,\tilde y_1)$. 
By the Hurwitz formula we obtain $g(\bar{D}) = 6$.
Now, since we are assuming ${\Xi^{0}}$ reducible, it follows that $D$, and then $\bar{D}$, is also reducible, so 
$\bar{D}= D_1 \cup D_2$ with $D_1 \cap D_2 = \{(\tilde y_1,\tilde y_2), (\tilde y_2,\tilde y_1)\}$ and each component is birational to $\tY$. Thus $D_i$ has 
genus $\geq g(\tY)=3$ for $i=1,2$, which implies that the arithmetic genus of $\bar{D}$ is at least 7. This gives a contradiction. 
Consequently, $\Xi$ is an irreducible genus 2 curve. 
\end{proof}


\begin{rem}
In \cite[Proposition 4.18]{lo} we saw how to find the curve $\Sigma$ with $P \simeq J\Sigma$ explicitely in term of the Weierstrass points of the curve $Y$. The same construction 
also works for covers of $R_2$ and $R_1$. Let $f: Y \ra X$ be a cover associated to an element of $R_2$ and let $\tilde f: \tilde Y \ra \tilde X$ be its normalization.
Moreover, let the notation be as in this section above. In particular $\tilde f(q_1)  = \tilde f(q_2) = \tilde f(q_3) = p_1$ and  $
\tilde f(q_4)  = \tilde f(q_5) = \tilde f(q_6) = p_2$. Considering the canonical map $\pi: \tY^{(2)} \ra \Pic^2(\tY)$, we can identify 
$\Pic^2(\tilde Y) \setminus H_{\tY}$ with $\tilde Y^{(2)} \setminus \pi^{-1}(H_{\tilde Y})$. Using this,
the Weierstrass points of $\Sigma$ are given by $[q_1 + q_2], [q_1 + q_3], [q_2 + q_3], [q_4 + q_5], [q_4 + q_6]$ and $[q_5 + q_6]$. An analogous construction works for
the elements of $R_1$.
\end{rem}

\section{The Prym map is finite}

We saw in \cite[Theorem 5.1]{lo} that the Prym map $Pr: {}_{S_3}\cM_2 \ra \cA_2$ is finite of degree 10 onto its image. Since ${}_{S_3}\cM_2$ is open dense in
${}_{S_3}{\widetilde \cM}_2$, the (extended) Prym map $Pr: {}_{S_3}{\widetilde \cM}_2 \ra \cA_2$ is also of degree 10. In this section we show that it does not 
have positive dimensional fibres in the boundary. 
We keep the notations of the previous section.

\begin{thm}
The Prym map $Pr: {}_{S_3} \widetilde{\cM}_2 \ra \cA_2$ is finite. 
\end{thm}

\begin{proof}
According to Theorem \ref{imagethm}, the Prym map is compatible with the coverings as given there. Hence, since $Pr: {}_{S_3} \widetilde{\cM}_2 \ra \cA_2$
is proper according to Theorem \ref{thm6.1}, it suffices to consider the 3 pieces separately and show that they have finite fibres. 
As mentioned above, this has been proved for $Pr: {}_{S_3}\cM_2 \ra \cJ_2^u$  already in \cite{lo}[Theorem 5.1].

The proof for $Pr: R \ra Pr(R) \subset \cJ^r_2$ is similar, but for sake of completeness we sketch it here.
Consider first the restriction of the Prym map to $R_2$.  Let $f:Y \ra X$ be a cover given by an element of $ R_2$ and
$Pr(f) = J\Sigma$ its image in $\cJ^r_2$. Let $\tY$ (of genus 3) and $\tX$ (elliptic curve) denote the normalizations. The corresponding map 
$\tilde{f} : \tY \ra \tX$ fits into a commutative diagram
\begin{equation} \label{commdiag}
\xymatrix{
\tY \ar[r]  \ar[d]_{\tilde f} & \PP^1 \ar[d]_{\bar{f}} & \Sigma \ar[l]_{\varphi} \ar[ld]^{\psi} \\
\tX  \ar[r]^{h_{\tX}} & \PP^1&
}
\end{equation}
where $\delta$ is the unique 2:1 map compatible with the hyperelliptic map of $\tY$.
The map $\bar{f}$ is doubly ramified at one point and simply ramified at 2 others. If $w_1,\ldots,w_6$ denote the Weierstrass points of the curve $\Sigma$, 
the 6:1 map $\psi = \varphi \circ \bar{f}$ (where $\varphi$ the hyperelliptic cover) is given by a pencil $g^1_6 \subset |3K_{\Sigma}|$, the ramification 
divisor of which consists of the 6 Weierstrass points of $\Sigma$, the 4 preimages of 2 ramification points over $h_{\tX}(p_3)$ and $h_{\tX}(p_4)$ and the 
2 preimages of the doubly ramified point over $h_{\tX}(\tilde x_1) = h_{\tX}(\tilde x_2)$. So one of the fibres of $\psi$ is $2w_i +2w_j +2w_k$ for some $1 \leq i < j <k \leq 6$ and
we denote the map $\psi$ by $\psi_{2(w_i + w_j + w_k)}$. Note that $\psi_{2(w_i + w_j + w_k)} = \psi_{2(w_l + w_m + w_n)}$ if $\{i,j,k,l,m,n\} = \{ 1, \dots,6\}$. 
In particular, $\psi$ has 5 distinct points in his branch locus.  

Conversely, the choice of 3 Weierstrass points $w_i,w_j,w_k$ such that $\bar f$ is doubly ramified at one point gives 5 marked points on $\PP^1$.
One can recover an element of $\Pr^{-1}(J\Sigma)$ as follows. Let $\bar y \in \PP^1$ the doubly ramified point of 
$\bar{f}$. Consider the uniquely determined elliptic curve $\tX$ given as a double cover $h_{\tX}: \tX \ra 
\PP^1$, branched over 4 of the marked points, all but $\bar x:=\bar{f}(\bar y)$. Define $\tilde f: \tY \ra \tX$ as the 
normalization of the pullback of $\bar{f}$ by $h_{\tX}$. 
One verifies that $\tY$ is a hyperelliptic curve (of genus 3) and $\bar{f}$ is non-cyclic and doubly ramified at the 2 points $\tilde y_1, \tilde y_2 $ over 
$\bar y= \bar{f}^{-1}(\bar x)$. 
Let 
$$
Y=\tY / \tilde y_1 \sim \tilde y_2 \quad \mbox{and} \quad  X=\tX / \tilde x_1 \sim \tilde x_2,
$$ 
where $h_{\tX}^{-1}(\bar x)= \{\tilde x_1, \tilde x_2\}$. Then the corresponding  
cover $f:Y \ra X$ defines an element of $R_2$ which maps to $J\Sigma$ under the Prym map. Since this construction depends only on the choice of the 3 Weierstrass points 
$w_i,w_j,w_k$, this implies that the map $Pr|_{R_2}: R_2 \ra \cJ^r_2$ has finite fibres. 

For the restriction of $Pr$ to $R_1$ the argument is similar and will be omitted (here the map $\bar f: \PP^1 \ra \PP^1$ has 2 doubly ramified points which lead to 
2 singular points of $Y$). Note that $R_1$ is connected of dimension 1, so it suffices to show that that
$Pr|_{R_1}$ is not constant. This completes the proof of the fact that the map $Pr: R \ra \cJ^r_2$ has finite fibres.  

It remains to show that the map $Pr: S \ra \cE_2$ is finite.
First let $f: Y \ra X$ be a cover associated to an element of $S_2$. 
Then $Y = Y^1 \cup Y^2$ (respectively $X = X^1 \cup X^2$) with smooth curves $Y^i$ of genus 2 (respectively $X^i$ of genus 1) intersecting 
transversally in 1 point and $f = f_1 \cup f_2$ with $f_i: Y^i \ra X^i$ doubly ramified at one point $y_i$ for $i = 1$ and 2.  
According to Proposition \ref{prop5.3} we have
$$
P(f) = P_1 \times P_2,
$$
where $P_i = \Ker (\Nm_{f_i}: JY_i \ra JX_i)$ for $i=1,2$.

It is enough to show that the Prym variety $P_i$ associated to the non-cyclic 3:1 
cover $Y^i \ra X^i$ defines a 1-dimensional family, when $X^i$ varies. Let $Z^i$ denote the Galois closure of $Y^i/X^i$.
The Galois group of $Z^i/X^i$ certainly is $S_3$.  Denoting by ''$\sim$'' isogeny and by $P(\cdot)$ the corresponding Prym varieties, we have
$$
JZ^i \sim P(Z^i/Y^i) \times P(Y^i/X^i) \times JX^i.
$$
On the other hand, according to \cite{rr},
$$
JZ^i \sim P(Y^i/X^i)^2 \times JX^i.
$$
This implies that it is enough to show that $P(Z^i/Y^i)$  defines a 1-dimensional family when $X^i$ varies.

Now $Z^i \ra Y^i$ is an \'etale double cover of a hyperelliptic curve
and for such a double cover it is very well known
that the Prym variety is a  product of 2 Jacobians, one of which may be 0 (see \cite{mu}). In our case,  $P(Z^i / Y^i) $ is isomorphic
to an elliptic curve $D_i$ which is an \'etale double cover of $X_i$. This proves the theorem for the 
elements in the image of $S_2$. The same argument works for the Prym varieties in the image of $S_1$,
since they are of the form $P=J\tY^1 \times P_2 $, where $P_2 = P(Y^2/X^2) $ (see Proposition \ref{prop5.4}).
We complete the proof by observing that $S_0$ is 0-dimensional. 
\end{proof}

\end{document}